\documentclass{gtpart}

\usepackage{stmaryrd}

\usepackage{anysize}
\marginsize{1in}{1in}{1in}{1in}

\usepackage{amsbsy}
\usepackage{amscd}
\usepackage[mathscr]{eucal}
\usepackage{verbatim}
\usepackage{version}
\usepackage{graphicx}
\usepackage{color}

\title{Refined curve counting on complex surfaces}

\author{Lothar G\"ottsche}
\givenname{Lothar}
\surname{G\"ottsche}
\address{International Centre for Theoretical Physics \\ Strada Costiera, 11 \\ Trieste, Italy}
\email{gottsche@ictp.it}
\urladdr{http://users.ictp.it/~gottsche}

\author{Vivek Shende}
\givenname{Vivek}
\surname{Shende}
\address{University of California, Berkeley \\ Department of Mathematics \\ 970 Evans Hall \\ Berkeley, California, USA 94720-3840}
\email{vivek@math.berkeley.edu}
\urladdr{http://math.berkeley.edu/~vivek}


\arxivreference{arXiv:1208.1973}
\arxivpassword{wq23m}

\volumenumber{}
\issuenumber{}
\publicationyear{}
\papernumber{}
\startpage{}
\endpage{}
\doi{}
\MR{}
\Zbl{}
\received{}
\revised{}
\accepted{}
\published{}
\publishedonline{}
\proposed{}
\seconded{}
\corresponding{}
\editor{}
\version{}

\definecolor{darkgreen}{rgb}{0.0, 0.7, 0.0}
\definecolor{darkpurple}{rgb}{0.7, 0.0, 0.7}

\def\Bbb{\mathbb}

\def\lbr{\llbracket}
\def\rbr{\rrbracket}

\newenvironment{pf*}[1]{\proof[#1]}{\endproof}

\newtheorem{thm}{Theorem}

\newtheorem{lem}[thm]{Lemma}
\newtheorem{cor}[thm]{Corollary}
\newtheorem{prop}[thm]{Proposition}
\newtheorem{conj}[thm]{Conjecture}

\theoremstyle{definition}
\newtheorem{quest}[thm]{Question}
\newtheorem{defn}[thm]{Definition}
\newtheorem{rem}[thm]{Remark}

\newtheorem{nota}[thm]{Notation}
\newtheorem{example}[thm]{Example}
\newtheorem{rec}[thm]{Recursion}

\newcommand{\lemref}[1]{Lemma~\ref{#1}}
\newcommand{\propref}[1]{Proposition~\ref{#1}}

\newcommand{\defref}[1]{Definition~\ref{#1}}
\newcommand{\remref}[1]{Remark~\ref{#1}}
\newcommand{\conjref}[1]{Conjecture~\ref{#1}}
\newcommand{\recref}[1]{Recursion~\ref{#1}}
\renewcommand{\C}{{\Bbb C}}
\renewcommand{\Z}{{\Bbb Z}}
\renewcommand{\Q}{{\Bbb Q}}
\renewcommand{\R}{{\Bbb R}}
\newcommand{\proj}{{\Bbb P}}
\newcommand{\CP}{\proj}

\newcommand{\LL}{\mathbb{L}}
\newcommand{\CC}{\mathbb{C}}

\newcommand{\ch}{\operatorname{ch}}

\newcommand{\res}{\mathop{\text{\rm res}}}

\newcommand{\Coeff}{\mathop{\text{\rm Coeff}}}
\def\oo{{\cal O}}

\def\P{{\mathbb P}}
\def\<{\langle}
\def\>{\rangle}

\newcommand{\Res}{\operatornamewithlimits{Res}}

\def\trop{{\text{trop}}}
\def\R{{\mathbb R}}

\def\Q{{\mathbb Q}}

\def\P{{\mathbb P}}
\def\Z{{\mathbb Z}}
\def\C{{\mathbb C}}

\def\oo{{\mathcal O}}

\def\={\;=\;}

\def\+{\,+\,} \def\-{\,-\,}
\def\be{\begin{equation}}
\def\ee{\end{equation}}
\def\bes{\begin{equation*}}
\def\ees{\end{equation*}}

\def\cC{\mathcal C}
\def\ZZ{\mathbb Z}

\newcommand{\que}{Q}
\newcommand{\gring}{\mathfrak{M}}
\newcommand{\mhm}{\mathrm{MHM}}

\begin{document}

\dedicatory{in memory of Friedrich Hirzebruch}


\begin{abstract}
	We define refined invariants which ``count'' nodal 
	curves in sufficiently ample linear systems
	on surfaces, conjecture that their generating function
	is multiplicative, and conjecture explicit formulas in the
	case of K3 and abelian surfaces.  We also give a refinement of
	the Caporaso-Harris recursion, and conjecture that it produces
	the same invariants in the sufficiently ample setting.  The
	refined recursion specializes at $y = -1$ to the Itenberg-Kharlamov-Shustin
	recursion for Welschinger invariants.  We find similar interactions
	between refined invariants of individual curves and real invariants
	of their versal families. 
\end{abstract}

\maketitle
\section{Introduction}

Given a general
elliptic fibration $K3 \to \mathbb{P}^1$,
we learn by computing $\chi(K3) = 24$ that there must be 24 nodal fibers.
For more general irreducible curve classes on a K3, 
Yau and Zaslow \cite{YZ} argued
that taking the Euler characteristic of the 
relative compactified Jacobian would again yield the number of maximally
degenerate fibers; the arguments were clarified by Beauville \cite{Bea}
and by Fantechi, van Straten, and the first author \cite{FGvS}.  

For more general families of curves, similar arguments may be made
in terms of the relative Hilbert schemes of points.  
Recall that for a smooth projective curve $C$ of genus $g$, 
the ``Macdonald formula'' asserts\footnote{We recall the derivation of this formula from \cite{M}. 
For any reasonable topological space we have 
$$H^*(X^{(n)})  = H^*(X^n/S_n) = H^*(X^n)^{S_n} = (H^*(X)^{\otimes n})^{S_n}.$$  
From this it follows that taking the generating function of Euler characteristics 
gives $$\sum_{n=0}^\infty q^n \chi (X^{(n)}) = \left( \sum_{n=0}^\infty q^n \right)^{\chi(X)}.$$
Taking $X$ to be a smooth curve $C$ and recalling that for these 
$C^{(n)} = C^{[n]}$ gives the formula stated.}
$$\sum_{n = 0}^\infty q^{n-g+1} \chi(C^{[n]})=\left(\frac{q}{(1-q)^2}\right)^{1-g}.$$

Let 
$\cC \to B$ be a family of reduced planar curves of arithmetic 
genus $g$, and let $\cC^{[n]} \to B$ be the relative Hilbert schemes.  
Certain string-theoretic ideas of Gopakumar, Katz, Klemm, and Vafa \cite{GV, KKV}
motivate the consideration of the 
following series, and the following change of variables:\footnote{
If $B$ has many connected components, one should perform the change of
variables component by component; on each $g$ should be intepreted as the arithmetic
genus of the fibre on the component in question.  One can avoid 
this unpleasantness by indexing as in \cite{PT2} and elsewhere
by the Euler characteristic of the ideal sheaf on the LHS and by the genus rather than cogenus
on the right; however, for our purposes indexing by the cogenus is far more convenient.}
\[\sum_{n=0}^\infty q^{n+1-g} \chi(\cC^{[n]}) = 
\sum_{i=0}^\infty n^i_{\cC/B} \cdot \left(\frac{q}{(1-q)^2}\right)^{i+1-g}.\]
If
in fact $B$ was a union of points and $\cC$ a disjoint union of 
smooth curves;
in this case the numbers $n_i$ would just count the number of curves 
of cogenus $i$.  In general we view the $n_i$ as the `virtual' number of curves 
of cogenus $i$ in the familiy $\cC \to B$.

This makes sense even when $\cC \to B$ is a single curve $C \to B = \mathrm{pt}$;
in this case we write simply $n^i_C$.  The Macdonald formula is equivalent
to the assertion that, when $C$ is smooth, $n^0_C = 1$ and $n^i_C = 0$ for all
$i > 0$. 
More generally, Pandharipande
and Thomas \cite{PT2} prove that $n^i_{C} = 0$ when $C$ has
cogenus $\delta(C) < i$.  Whenever the singularities of $C$
are unions of smooth branches, the last nonvanishing term is
calculated either from \cite{Bea} or \cite{FGvS} to be $n^i_C = 1$. 

In the relative situation, it is helpful to view $n^i$ as a constructible function 
$n^i:B \to \ZZ$ given by $b \mapsto n^i_{\cC_b}$; evidently 
$n^i_{\cC_B/B} = \int_B n^i(b) d \chi(b)$.  

In addition to the naive or virtual interpretation above, in good cases
the $n^i$ carry actual enumerative meaning.  Suppose that $i$ is the
maximum cogenus of any curve in the family,
that there are finitely many curves of cogenus $i$, 
and that all these curves are immersed.  Then
$n^i_{\cC/B}$ is just the number of these curves. 
One exploits this observation by finding another
way to express the Euler characteristics of the relative Hilbert schemes.  
In particular,
the following two results have recently been established: 

\begin{thm} \label{thm:mult}  \cite{Sh}.  Fix $C$ a reduced plane curve,
  $\Lambda$ a versal deformation of its singularities, and 
  $\Lambda^\delta \subset \Lambda$ the locus of cogenus $\delta$ curves. 
  Then $n^\delta_C$ is the multiplicity of $\Lambda^\delta$. 
\end{thm}

\begin{thm} \label{thm:univ} \cite{KST}.  
  Let $S$ be a surface, and $L$ a $\delta$-very-ample line bundle.  
  Then the number of $\delta$-nodal curves in a general
  $\proj^\delta \subset |L|$ is $n^\delta_{\cC/\proj^\delta}$, which
  moreover
  is given by a certain
  explicit combination of integrals of Chern classes
  of $L^{[n]}$ and $TS^{[n]}$. 
\end{thm}

The significance of the second result is that, due to a theorem of
Ellingsrud, Lehn, and the first author \cite{EGL}, such integrals depend
in a universal way on the Chern classes of $L$ and $S$.  Such universality
of the counts of nodal curves had previously been conjectured
by the first author \cite{G}, and previously proven by \cite{Tz}.
Note that in particular, this universality implies that the numbers
$n^\delta_{\cC/\proj^\delta}$ will not vary under a deformation of the pair
$(S, L)$ which preserves the ampleness condition, since the integrals
of Chern classes of $S, L$ will not change.  Likewise, Theorem \ref{thm:mult} 
plus the fact that Euler numbers of the Hilbert schemes of singular curves 
may be recovered from the HOMFLY polynomials
of their links \cite{Mau} implies that the numbers $n^\delta_C$ do not 
vary under an equisingular deformation. 

\vspace{4mm}

The present article poses the following question: {\em does
replacing the topological Euler characteristic on the left
hand side by more sophisticated invariants have an enumerative
counterpart on the right?}  

\vspace{4mm}

We begin in Section \ref{sec:curve} by studying the case of a single curve.  
At the outset, we work in the Grothendieck ring of varieties, and consider
$\sum_{n=0}^\infty q^{n+1-g} [C^{[n]}] $.  However, this incorporates 
global information we would rather not consider, and in particular depends
on the motive of the normalization $\widetilde{C}$.  We can remove 
the global contributions by dividing out by the analogous series for
$\widetilde{C}$.  The quotient only depends on the singularities.
We write $\widetilde{g} = g - \delta$ for the genus of $\widetilde{C}$, and show
there exist 
classes $\widetilde{N}^i_C$ in the Grothendieck ring of varieties such that

\[ 
\frac{\sum_{n=0}^\infty q^{n+1-g} [C^{[n]}] }{\sum_{n=0}^\infty q^{n+1-\widetilde{g}} [\widetilde{C}^{[n]}]}
= \sum_{i=0}^\delta \widetilde{N}^i_C \cdot \left( \frac{q}{(1-q)(1-q [\mathbb{A}^1] )} \right)^{i-\delta}
\] 

The right hand side moreover splits into a product over the singularities of $C$. 
See Equation \eqref{eq:product} and Corollary \ref{cor:vanishing} below for proofs and further discussion.

According to Theorem \ref{thm:mult}, the Euler number $\chi(\widetilde{N}^i_C)$ gives the 
multiplicity of a certain locus, and is in particular positive.  In examples we see a stronger
positivity: 

\begin{conj}
	$\widetilde{N}^i_C \in \ZZ_{\ge 0}[\mathbb{A}^1]$.
\end{conj}

This conjecture is verified computationally for singularities
of the form $x^p = y^q$ where $(p,q) = 1$ and $p < 12, q < 20$ using
the formulas of \cite{ORSh} for the classes of the Hilbert schemes.  

The meanings of the $\widetilde{N}^i_C$ remain mysterious, but there is some
evidence that they may be related to geometry over $\R$.  We define real analogues
$n^{i,\R}$ of the $n^i$ by using the compactly supported Euler number of the real
locus, and show these again vanish for $i > \delta(C)$.  
These have an interpretation analogous to Theorem \ref{thm:mult}. 
Recall that nodes of real curves
come in three types: elliptic ($x^2 + y^2 = 0$), hyperbolic ($x^2 - y^2 = 0$), 
and pairs of complex conjugate nodes. 
Thus in the real deformation, the loci $B_+^k$ of $k$-nodal curves split into 
components according to the types of the nodes.  

\begin{thm}
  Let $C$ be a real reduced plane curve, and let $\cC \to B$ be a locally 
  versal deformation of its singularities.  Let $B^{\delta, \delta_-}_+$
  be the locus of nodal curves with $\delta$ nodes of which $\delta_-$ 
  are hyperbolic.  Let $D^j$ be a general real disc of dimension $j$
  passing near $[C] \in B$. 
  Then \[n^{j,\R}_C = \sum_i (-1)^i D^j \cap B^{j, i}.\]
\end{thm}

Note in particular that while the individual terms on the RHS of the formula may depend
on the location of the disc, the theorem asserts that their sum does not.  

For the simple singularities, we give a combinatorial formula
for the $\widetilde{N}^i_C$ in terms of 
the Dynkin diagram.  The formula refines the analogous prescription
for the multiplicities given in \cite{Sh}.  Geometrically, this may be interpreted as
the choice of a particular real form (in the unibranch case, there is no choice) and a particular 
disc $D$ in the above statement so that the coefficient
of $\LL^i$ in $N^j$ is $D^j \cap B^{j,i}$.  For $j = \delta$, Duco van Straten has
conjectured that such a disc may be found for any singularity. 

\vspace{2mm} 
In Section \ref{sec:families}, we turn to the case of linear systems of curves on
surfaces.  From the point of view of the argument in \cite{KST}, it is natural
to refine the Euler characteristic to Hirzebruch's $\chi_{-y}$ genus,
 since the
latter both factors through the Grothendieck ring of varieties and may be calculated
in terms of Chern classes. 
Recall that the $\chi_{-y}$ genus is given by $\chi_{-y}(X)=\sum_{p,q} (-1)^{p+q} {\mathrm h}^{p,q}(X)$ where the ${\mathrm h}^{p,q}(X)$ are the Hodge numbers.
 We define invariants 
$N^i_{\cC/\proj^\delta} \in \ZZ[y]$ which refine the $n^i_{\cC/\proj^\delta}$ 
of Theorem \ref{thm:univ} above.  As before, we are motivated by the MacDonald formula,
which in this case reads 
$$\sum_{n = 0}^\infty q^{n-g+1} \chi_{-y}(C^{[n]})=\left(\frac{q}{(1-q)(1-qy)}\right)^{1-g}.$$
Thus we {\em define} the $N^i_{\cC/B}$ by equating terms in Laurent
series expansions of the following (see \defref{def:change}):

\[\sum_{n=0}^\infty q^{n+1-g} \chi_{-y}(\cC^{[n]}) = 
\sum_{i=0}^\infty N^i_{\cC/B} \cdot \left(\frac{q}{(1-q)(1-qy)}\right)^{i+1-g}.\]

As before, if $B$ were just isolated points and $\cC$ were a collection
of smooth curves, then $N^i_{\cC/B}$ would be just the number of 
curves of cogenus $i$.
In particular, for a single smooth curve, we have $N^0_C = 1$ and
$N^i_C = 0$ for $i > 0$.  More generally, the $N^i_C$ vanish for
$i > \delta(C)$; this is particular to the $\chi_{-y}$ genus: the analogous statement does not hold for the virtual Poincar\'e 
polynomials.\footnote{This does not contradict the statement above that the $\widetilde{N}^i_C$ vanish even
in the Grothendieck ring of varieties; the point being that the $N^i_C$ and $\widetilde{N}^i_C$ 
have different virtual Poincar\'e polynomials but the same $\chi_{-y}$ genus.  Ultimately 
this is because the $\chi_{-y}$ genus of an abelian variety vanishes.}

Let $S$ be a surface, $L$ a line bundle on it, $\CP^\delta \subset |L|$ 
a linear system. 
Let $S^{[n]}$ be the Hilbert scheme of $n$ points on $S$, 
and let $Z_n(S)\subset S\times S^{[n]}$ be the universal family, 
with the projections $q:Z_n(S)\to S$, $p:Z_n(S)\to S^{[n]}$. 
Let $L^{[n]}:=p_*q^*L$.  This is a vector bundle of rank 
$n$ on $S^{[n]}$ with fibre $\mathrm{H}^0(Z,L|_Z)$ over $Z\in S^{[n]}$. 
Let $\cC$ be the universal curve over $\CP^\delta$ 
and  $\cC^{[n]}$ 
the relative Hilbert scheme of points.  
It  is the scheme theoretic zero locus of a tautological section 
of $L^{[n]}\boxtimes H$; when $\cC^{[n]}$ 
is nonsingular this section is transverse. This allows to compute $\chi_{-y}(\cC^{[n]})$ as an intersection number on $S^{[n]}$
(see \propref{prop:hilbgenus}).

Experimental evidence suggests: 

\begin{conj} \label{conj:vanishing}
	Let $L$ be a line bundle on a surface $S$, and $\proj^\delta \subset |L|$ a linear subsystem
	with tautological curve $\cC \to \proj^\delta$.   Assume that $\proj^\delta$ contains
	no non-reduced curves, and that the total space of the relative
	Hilbert scheme $\cC^{[n]}$ is smooth for all $n$.  Then $N^i_{\cC/\proj^\delta}$ vanishes for
	$i > \delta$.  
\end{conj}

Unlike in the Euler characteristic setting, one cannot prove this by ``integrating
over the base.'' Indeed, already in the smooth case ($\delta = 0$)  Pandharipande and Fantechi \cite{Pan} 
have found families of smooth curves over
a smooth base curve with nonzero invariants $N^i$ for some $i > 0$.  
The base curves 
are always of genus $> 0$, and indeed no such examples can exist over a
simply connected base.  We are not sure what further aspects of the geometry of families
of curves on surfaces are implicated.  A calculation of Migliorini has shown that the vanishing
cannot be expected for the analogous expression involving Hodge polynomials, already over a one
dimensional base.  Nevertheless, we have been able to show

\begin{thm} \label{thm:introvanishing}
	The conjecture holds when $K_S$ is numerically trivial. 
\end{thm}

We also give some evidence for the general statement.  By \cite{EGL}, one can reduce the validity of the conjecture
to the case of $\proj^2$ and $\proj^1 \times \proj^1$.  Here we may calculate
in low degrees by equivariant localization.

In Section \ref{sec:deepest}, we focus on the invariants $N^\delta_{\cC/\proj^\delta}$.
Assuming Conjecture \ref{conj:vanishing}, this is the last nonvanishing $N^i$. Under this assumption we show  that there are universal power series $A_1,\ldots, A_4$  in $\mathbb{Q}[y]\lbr s\rbr $, such whenever $L$ is a $k$-very ample line bundle on a surface $S$  we have
\[\sum_{\delta = 0}^\infty N^\delta_{\cC/\proj^\delta} s^\delta \equiv A_1^{L^2} A_2^{L K_S} A_3^{K_S^2} A_4^{c_2(S)}
+ O(s^{k+1}).\]
We define polynomials $N^\delta_{\delta,[S,L]}$ by 
\[\sum_{\delta = 0}^\infty N^\delta_{\delta,[S,L]} s^\delta=A_1^{L^2} A_2^{L K_S} A_3^{K_S^2} A_4^{c_2(S)}.\]
In particular $N^\delta_{\delta,[S,L]}=N^\delta_{\cC/\proj^\delta}$, if $L$ is $\delta$-very ample.
The brackets in the notation 
serve to remind us that it depends only on the cobordism class of $(S,L)$.

As in \cite{G}, it is easiest to express the $A_i$ after a change of variable.
Consider the following series in $\Q[y,y^{-1}]\lbr q\rbr $
\begin{align*}
\widetilde\Delta(y,q)&:= q \prod_{n=1}^{\infty}(1-q^n)^{20}(1-yq^n)^{2}
(1-y^{-1}q^n)^2,\\
\widetilde{DG}_2&:=\sum_{m= 1}^\infty m q^{m} \sum_{d | m} \frac{[d]_y^2}{d} 
\end{align*}
and let $D=q\frac{d}{dq}$.  Above $[n]_y:= \frac{y^{n/2}- y^{-n/2}}{y^{1/2}-y^{-1/2}}$. 

It is also convenient to introduce the notation $\overline{N}^\delta_{\delta,[S,L]}
:= y^{-\delta} N^\delta_{\delta,[S,L]}$; these invariants are symmetric under $y \to 1/y$.  

\begin{conj} \label{conj:Bs}
There exist power series 
$B_1(y,q)$, $B_2(y,q)$ in $\Q[y,y^{-1}]\lbr q\rbr $, such that 
\[
\sum_{\delta\ge 0}
\overline{N}^\delta_{\delta,[S,L]}(y)\widetilde{DG}_2^\delta=\frac{(\widetilde{DG}_2/q)^{\chi(L)}B_1(y,q)^{K_S^2}B_2(y,q)^{LK_S}}
{\big(\widetilde\Delta(y,q)\, D \widetilde{DG}_2/q^2)^{\chi(\oo_S)/2}}.
\]
\end{conj} 
The first 11 terms of the series $B_1, B_2$ are given explicitly in Section \ref{sec:deepest}.

Assuming Conjecture \ref{conj:vanishing}, the content of the above assertion is exhausted 
by the case where $S$ is a K3 surface.  Moreover as in \cite{G} it may be reformulated without
the expansion in powers of $\widetilde{DG}_2$:

\begin{conj}\label{conj:K3A} Let $(S_g, L_g)$ be K3 surfaces of genus $g$ with irreducible polarizations. 
Then for any $k$, 
\[ \sum_{g=k}^\infty q^{g-1} \overline{N}^{g-k}_{g-k,[S_g,L_g]}(y)  = 
\frac{\widetilde{DG}_2(y,q)^k}{\widetilde{\Delta}(y,q)}.\]
Similarly, if $(A_g, L_g)$ are abelian surfaces of genus $g$ with irreducible polarizations, 
\[
\sum_{g=k+2}^\infty \overline{N}^{g-k-2}_{g-k-2,[A_g,L_g]} q^{g-1}  =  \widetilde{DG}_2(y,q)^k D\widetilde{DG}_2
(y,q).
\] 
\end{conj} 

\begin{rem}
After the original version of this paper, we have 
proven  \conjref{conj:K3A},  \conjref{conj:K3H}, and given more general formulas in  \cite{GSh}. 
In the current paper we show that 
\conjref{conj:K3A} follows from its validity for {\it both} K3 and abelian surfaces at $k=0$, and show the $k=0$ case for K3 surfaces. 
The proof uses the existence of K3 and abelian surfaces of all genera, and the multiplicative nature of
the formulas. The proof of the conjecture in \cite{GSh} relies on these results, and consists in 
establishing the validity of the conjecture at $k=0$ for abelian surfaces by studying moduli spaces of pairs on abelian surfaces
and a crucial use of Theorem \ref{thm:introvanishing}.  
\end{rem}

Conjecture \ref{conj:K3A} would also follow from its validity for all $k$ for the K3 surface alone. 
Here we note a remarkable coincidence: the series on the RHS of the above formula
appears in the work of \cite{MPT} on computing descendant invariants in the (reduced)
Gromov-Witten or stable pairs theory of a K3 surface.  This leads to a further reformulation: 


\begin{conj} \label{conj:K3H}
  Let $(S, L)$ be a irreducibly polarized $K3$ surface of genus $g$, and let $H$ be the hyperplane
  class on $|L|$.  Then for all $k$, 
  \[ 
  (y - 2 + y^{-1})^{k-1} \overline{N}^{g-k}_{g-k, [S,L]} = \sum_{n=0}^\infty 
  y^{n+1-g} \int_{\cC^{[n]}_{|L|}} c_{n+g-k}(T \cC^{[n]}_{|L|})\cdot \rho^*(H^k).
  \]
\end{conj}

\vspace{2mm}
Thus far we have been discussing curves with a number of nodes small compared to the ampleness of
the line bundle
$L$; this is the regime to which the conjectures of \cite{G} and the arguments of \cite{KST} apply.
However, when the surface is $\proj^2$, the recursion of Caporaso and Harris
\cite{CH} determines the degrees of {\em all} such loci of nodal curves, without any such restriction on 
the ampleness.  Indeed, it determines more: fix a line $H \subset \proj^2$, and sequences $\alpha, \beta$
of integers specifying respectively fixed and moving tangency conditions to $H$.  
Then the Caporaso-Harris
recursion determines the degrees $n^{d,\delta}(\alpha,\beta)$ of the loci of curves with $\delta$
nodes and satisfying the tangency conditions $\alpha, \beta$.\footnote{Vakil has 
generalized the Caporaso-Harris recursion to the case of rational ruled surfaces.
In Section \ref{sec:Severi} we treat these as well; we have restricted in the introduction to $\proj^2$ just
for ease of notation.}

In Section \ref{sec:Severi} we study the following formal refinement of the Caporaso-Harris recursion. 
We take from \cite{CH} the notation $I \alpha = \sum i \alpha_i$ and $|\alpha| = \sum \alpha_i$; note that
the curves counted by $n^{d,\delta}(\alpha, \beta)$ have degree $I\alpha + I\beta$. 

\begin{defn}\label{Nddef}
	The polynomials $\overline{N}^{d,\delta}(\alpha, \beta) \in \ZZ[y^{1/2},y^{-1/2}]$ are defined by the following recursion: 
\[
	\overline N^{d,\delta}(\alpha,\beta)=\sum_{k:\beta_k>0} [k]_y 
	\cdot \overline N^{d,\delta}(\alpha+e_k,\beta-e_k) +\sum_{\alpha',\beta',\delta'}
	\left( \prod_i [i]_y^{\beta_i'-\beta_i} \right) \binom{\alpha}{\alpha'}\binom{\beta'}{\beta} \overline N^{d-1,\delta'}(\alpha',\beta').
\]
\end{defn}

The limits on the sum and the initial conditions are the same as for the Caporaso-Harris recursion
and are given explicitly in Section \ref{sec:Severi}.  
The refined recursion immediately specializes to the Caporaso-Harris recursion upon setting $y = 1$, so certainly
$n^{d,\delta}(\alpha, \beta) = \overline{N}^{d,\delta}(\alpha,\beta)|_{y = 1}$.   On the other hand, we know that for
$d \gg \delta$ the Severi degrees $n^{d,\delta}((0,0,\ldots),(d,0,\ldots))$ are given
by the universal formulas, i.e., equal to the numbers $n^{\delta}_{\delta,[\proj^2, \mathcal{O}(d)]}$.  
We conjecture a refined analogue: 

\begin{conj} \label{conj:CHrefined}
   For $\delta \le 2d-2$,  we have 
   $\overline{N}^{d,\delta}((0,0,\ldots),(d,0,\ldots)) = y^{-\delta} N^\delta_{\delta,[\proj^2, \mathcal{O}(d)]}$.
\end{conj}

The equality at $y = 1$ follows from \cite{KS}. 
At $y=0$, the recursion
simplifies, allowing the RHS to be calculated explicitly; on the other hand, a result of Scala 
allows the LHS to be calculated as well \cite{Sc}; the answers match.  
Empirically we have verified the equality for some small $d, \delta$.  

\vspace{2mm}

In Section \ref{sec:RRT}, we note a connection to real enumerative geometry and to some ideas from tropical geometry. On a real toric surface $S$, there are real enumerative invariants, counting real $\delta$-nodal curves with suitable signs, the real analogues of the Severi degrees.
If $S$ is an unnodal del Pezzo surface, they coincide with the Welschinger invariants,  real analogues of the Gromov-Witten invariants. 
Mikhalkin \cite{Mi} has shown that the Severi degrees and the real enumerative invariants can be computed via tropical geometry: he introduces tropical Severi degrees and tropical Welschinger invariants by assigning Gromov-Witten and Welschinger multiplicities to tropical curves, and shows that they coincide with the Severi degrees and the real enumerative invariants respectively. 
The Caporaso-Harris formula has been derived tropically by Gathmann and Markwig \cite{GM}, and an analogue for the tropical Welschinger invariants by Itenberg, Kharlamov, and Shustin
\cite{IKS}.   These  are  specializations of the above recursion, specialized
at $1$ and $-1$ respectively. In particular the refined Severi degrees specialize to the Severi degrees and the tropical Welschinger invariants. 
In \cite{BG}, Block and the first author define and study tropical refined Severi degrees by assigning polynomial multiplicities to the tropical curves which specialize to the Gromov-Witten and Welschinger multiplicities.

On the right hand side of Conjecture \ref{conj:CHrefined}, the specialization $y \mapsto -1$ has an entirely
different meaning: we are taking the signatures of the relative Hilbert schemes and rearranging them in a certain way. 
On the other hand, a signed count of real nodal curves in a general $\proj^\delta$ is obtained from the 
$n^{\delta, \R}_{\proj^\delta}$ defined earlier.  In order that $n^{\delta, \R}_{\proj^\delta}$ match 
$\overline{N}^\delta_{\proj^\delta}(-1)$, the following property would suffice:

\begin{conj}
Let $\proj^\delta \subset |\mathcal{O}_{\proj^2}(d)|$ be a linear system determined by a subtropical collection
of real points. If $\cC^{[n]}/\proj^\delta$ is smooth, then its signature is equal to the Euler
characteristic of its real locus. 
\end{conj}
Here roughly speaking a collection of points in $(\R^*)^2$ is called subtropical, if it can be degenerated to a tropical collection of points without crossing walls, for the precise definition see \cite[Lemma 2.7.(3)]{IM}.

\vspace{4mm}

We remark briefly on related work.  In the physics literature there is the notion of the
refined topological string, which gives in some cases a one-parameter deformation of 
the various curve counting invariants \cite{IKV}.  Notably it does not have a ``worldsheet'' definition,
even in the sense of physics.  Mathematically, the refined theory is supposed to correspond \cite{DG} to 
the motivic DT theory \cite{KoS}; the lack of a worldsheet definition corresponds to the fact that 
we do not know how to correspondingly refine the Gromov-Witten invariants.  (For further
discussions of motivic DT theory, see, e.g. \cite{KoS2, BBS, CKK, BJM}.) There have also
been intimations that a specialization of the refined theory is related to real invariants \cite[Sec. 5]{KW}.  Our 
approach falls roughly into this paradigm insofar as $\overline{N}^i_{\cC/\proj^\delta}$ are assembled
from $\chi_{-y}$ genera of relative Hilbert schemes, which under the relevant assumptions on
$\proj^\delta$ are just the same as stable pairs spaces.  It might plausibly be hoped that the refined Severi degrees
also admit an interpretation in the stable pairs theory \cite{PT1, PT2} or its surface variant \cite{KT}.

\vspace{4mm}
\noindent {\bf Acknowledgements}.  We thank 
Paolo Aluffi, Florian Block, Barbara Fantechi, 
Ilia Itenberg, Davesh Maulik, Luca Migliorini, Grigory Mikhalkin, 
Hiraku Nakajima, Andrei Okounkov, Rahul Pandharipande, Sam Payne, 
Duco van Straten, Richard Thomas, Kota Yoshioka, 
and Don Zagier for helpful discussions.  The first named author thanks Friedrich Hirzebruch for interesting and useful discussions on the subject of this paper, and  also for his help and guidance over many years.

\vspace{4mm} 

\noindent {\bf Notation.}

We denote quantum numbers as 
\[ [n]_y := \frac{y^{n/2}- y^{-n/2}}{y^{1/2}-y^{-1/2}}. \]

By the Hirzebruch genus $X_{-y}$ we mean the characteristic class 
which on a bundle $E$
with Chern roots $x_i$ takes the value 
\[X_{-y}(E) = \prod \frac{x_i(1-ye^{-x_i(1-y)})}{(1-e^{-x_i(1-y)})} \in 
1+(x_i)\Q[y]\lbr x_i\rbr . \] 
Also let 
\[ch_{-y}(E) = \sum e^{x_i(1-y)}.\]
Setting by definition
\[\chi_{-y}(X,E):=\sum (-y)^p  \chi(X, \Omega^p \otimes E)  \]
we have according to Hirzebruch
\[\chi_{-y}(X,E) = \int_X ch_{-y}(E) X_{-y}(TX).\] 

When $E = \oo_X$ we suppress it from the notation.  Note that $\chi_{-y}(X)=\sum_{p,q}(-1)^{p+q}y^q 
\mathrm{h}^{p,q}(X)$ where $\mathrm{h}^{p,q}(X)$ are the Hodge numbers of $X$. Note the specializations
to topological Euler characteristic $\chi_{-1}(X) = \chi(X)$, holomorphic Euler characteristic 
$\chi_0(X) = \chi(X, \oo_X)$, and signature 
$\chi_1(X) = \sigma(X)$. 

\vspace{4mm} 

\noindent{\bf Glossary of notations}

\begin{tabular}{ll}
$X^{(n)}$ &$n$-th summetric power\\
$X^{[n]}$ &Hilbert scheme of $n$ points\\
$\cC^{[n]}$ &relative Hilbert scheme of points of family of curves\\
$\P^\delta\subset |L|$ &$\delta$ dimensional sublinear system of complete linear system\\
$\gring$ & Grothendieck ring of varieties\\
$[X]$ &class in $\gring$\\
$[S,L]$&cobordism class of a surface with a line bundle\\
$D^{S,L}(y,x,q)$ &generating function for integrals on $S^{[n]}$ which determines the $N^i_{\delta,[S,L]}$  \eqref{DSL}\\
$\que$ & $q/((1-q)(1-yq))$\\
$m(C)$ &tropical multiplicity of a tropical curve\\
$r(C)$& $[m(C)]_{-1}$ Welschinger multiplicity of a tropical curve\\
$M(C)$ &$[m(C)]_y$ refined multiplicity of a tropical curve\\
\end{tabular}

\smallskip
{\it Curve counting invariants defined via generating functions of invariants of (relative) Hilbert schemes of points:}\\
\begin{tabular}{ll}$n_C^i$ &invariants from Euler numbers for single curve  \\
$n^i_{\cC}$ & the same for family of curves\\
$n_C^{i,\R}$ & real invariant, from Euler number of real locus for single curve\\
$n_{\cC}^{i,\R}$ &the same for family of real curves\\
$\widetilde N^{i}_C$& invariants in Grothendieck group of varieties for single curve (see \eqref{eq:product})\\
${\mathcal N}^i_{\cC/B}$ &invariants for family of curves in mixed Hodge modules \eqref{eq:nyMHM}\\
${\mathbf n}^i_{\cC/B}$ &invariants for family of curves in Chern-Schwarz-Macpherson classes\\
${\mathbf N}^i_{\cC/B}$ &Brasselet, Sch\"urmann, Yokura refinement of the above\\
$N_\cC^i$, $N_{\cC/B}^i$  &invariants for family of curves from $\chi_{-y}$-genus  (see \defref{def:change})\\
$N_{\cC/\P^\delta}$ &the above for linear system $\P^\delta\subset |L|$\\
$N^i_{\delta,[S,L]}$ &polynomial in $L^2$, $LK_S$, $K_S^2$, $c_2(S)$ that
equals $N_{\cC/\P^\delta}$ for $L$ sufficiently \\ &ample (see \defref{Ndef})\\
${\mathbf N}^i_{\delta,[S,L]}$ & the same for $N_{\cC/\P^\delta}$ (see \defref{Ndef})\\
$\overline N^i_{\delta,[S,L]}$ &$N^i_{\delta,[S,L]}/y^{\delta}$\\
\end{tabular}
\smallskip

{\it Generalized Severi degrees:}\\
\begin{tabular}{ll}
$n^{d,\delta}$ &Severi degree on $\P^2$: number of $\delta$-nodal curves in $\P^2$ of degree $d$\\
$n^{d,\delta}(\alpha,\beta)$ &relative Severi degree:  curves with contact conditions\\
$n^{L,\delta}$ &Severi degree on toric surface \eqref{Srec}\\
$n^{L,\delta}(\alpha,\beta)$ &relative Severi degree \eqref{Srec}\\
$N^{d,\delta}$, $N^{d,\delta}(\alpha,\beta)$ & (relative) refined Severi degree on $\P^2$ \\
$N^{L,\delta}$,  $N^{L,\delta}(\alpha,\beta)$& (relative) refined Severi degree on toric surface \eqref{refCHrec}\\
$\overline N^{d,\delta}$, $\overline N^{d,\delta}(\alpha,\beta)$&  normalized refined Severi degrees (see \defref{Nddef})\\
$\overline N^{L,\delta}$,  $\overline N^{L,\delta}(\alpha,\beta)$& (see \defref{barNL})\\
$N^{L,\delta}_0$&  refined Severi degree on $\P^2$ counting irreducible curves\\
$W^{L,\delta}$& real  curve counting invariants of toric surface\\
$W^{L,\delta}_0$& real curve counting invariant for irreducible curves\\
$W^{d,\delta}_{\trop}$& tropical Welschinger invariants of degree $d$ on $\P^2$\\
$W^{L,\delta}_{\trop}$& tropical Welschinger invariants on toric surface\\ 
$\overline N^{L,\delta}_{\trop}$& refined tropical Severi degrees on toric surface\\ 
\end{tabular}

\section{Invariants of a single curve} \label{sec:curve}

\subsection{Refined invariants}

Let $\gring$ denote the Grothendieck ring of varieties; let $\LL = [\mathbb{A}^1]$ denote the class of the affine line.  
Kapranov \cite{Kap} introduced
the motivic zeta function of a variety: 
$$\zeta_X(q) = \sum_{n=0}^\infty \mathrm{Sym}^n X \cdot q^n \in \gring\lbr q\rbr.$$
When $X$ is a smooth proper curve of genus $g$, he showed that in fact
$(1-q)(1-q \LL)\zeta_X(q)$ is a polynomial of degree $2g$, and that moreover
one has a functional equation 
$$\zeta_X(1/q\LL) = \LL^{1-g} q^{2-2g} \zeta_X(q).$$

Motivated by a circle of ideas relating curve counting, Hilbert schemes on 
singular curves, and knot invariants \cite{GV, KKV, PT1, PT2, Sh, KST, OSh, ORSh, MY, MS, Mau}, we consider
a ``zeta function'' defined using the Hilbert schemes rather than the symmetric products.  The analogous
rationality and functional equation continue to hold.\footnote{This was observed for the Euler characteristics
in \cite[Prop. 3.13]{PT2}; the argument however is identical to that for the smooth case as in \cite{Kap} which
in turn is essentially the same argument as in Schmidt's original proof of the rationality and functional
equation for the zeta function of a curve \cite{Sch}.  We nevertheless include (again the same) proof for completeness.}


\begin{lem} \label{lem:hart}  \cite{H}.
  Let $C$ be a Gorenstein curve, and let $F$ 
  be a torsion free sheaf on $C$. 
  Write $F^*$ for $\mathcal{H} om (F,\mathcal{O}_C)$. 
  Then $\mathcal{E}xt^{\ge 1}(F,\mathcal{O}_C) = 0$ and 
  $F = (F^*)^*$.  Serre duality
  holds in the form $\mathrm{H}^i(C, F) = \mathrm{H}^{1-i}(C, F^* \otimes \omega_C)^*$.
  For $F$ rank one and torsion free, define its 
  degree $d(F):=\chi(F) - \chi(\mathcal{O}_C)$.   This satisfies 
  $d(F) = - d(F^*)$, and, for $L$ any line bundle, 
  $d(F \otimes L) = d(F) + d(L)$. 
\end{lem}

\begin{prop}  
\label{prop:ratfe}
Let $C$ be a complete, reduced,  irreducible Gorenstein, complex curve of arithmetic genus $g$. 
Then $f_C(q) = (1-q)(1-q\LL) \sum_d q^d [C^{[d]}]$ is a polynomial of degree
$2g$, satisfying $q^{2g} \LL^g f_C(1/q\LL) = f_C(q)$. 
\end{prop}
\begin{proof}
  Fix a degree one line bundle $\mathcal{O}(1)$ on $C$. 
  Let $\overline{J}^0(C)$ denote the moduli space of rank one, degree zero, torsion free sheaves \cite{AK}. 
  We map $C^{[d]} \to \overline{J}^0(C)$ by 
  sending the ideal $I \subset \mathcal{O}_C$ to the
  sheaf $I^* = \mathcal{H} om(I,\mathcal{O}_C) 
  \otimes \mathcal{O}(-d)$; the fibre is 
  $\proj(\mathrm{H}^0(C,I^*))$.  For $F \in \overline{J}^0(C)$, we write the Hilbert function as 
  $h_F(d) := \dim \mathrm{H}^0(C,F \otimes \mathcal{O}(d))$.  
   
  Fix $h = h_F$ for some $F$.  Evidently $h$ 
  is supported in $[0,\infty)$, and by Riemann-Roch and
  Serre duality is equal to
  $d+1-g$ in $(2g-2,\infty)$.  Inside $[0,2g-2]$, it
  either increases by 0 or 1 at each step.  Let $\phi_\pm(h)
  = \{ d\,|\, 2 h(d-1) - h(d-2) - h(d) = \pm 1\}$; evidently
  $\phi_- \subset [0,2g]$ and $\phi_+ \subset [1,2g-1]$.  Consider
  \begin{align*}
  j_h(q)&:=
  (1-q) (1-q\LL) \sum_{d=0}^\infty q^d [\proj^{h(d)-1}] = \sum_{d=0}^\infty q^d \left([\proj^{h(d)-1}]-(1+\LL)[\proj^{h(d-1)-1}]+\LL[\proj^{h(d-2)-1}]\right)\\&=
  \sum_{d \in \phi_-(h)} q^d \LL^{h(d)-1} -
  \sum_{d \in \phi_+(h)} 
  q^{d} \LL^{h(d-1)}. 
   \end{align*}
  This is a polynomial in $q$ of degree at most $2g$. 
  
  Now let $G = F^* \otimes \omega_C \otimes \mathcal{O}(2-2g)$, 
  and $h^\vee = h_G$. 
  By Serre duality and Riemann-Roch, $h^\vee(d) = h(2g-2-d) + d + 1 - g$, so
  in particular, $d \in \phi_{\pm}(h^\vee) \iff 2g-d \in \phi_{\pm}(h)$.  
  It follows that $q^{2g} \LL^g j_{h^\vee}(1/q\LL) = j_h(q)$. 
  
  Finally,
  stratify 
  $\overline{J}^0(C)$ into strata over which $h_F$ is constant.  The restriction of 
  $C^{[d]}$ to each stratum is the projectivization of a vector bundle of rank $h_F(d)$.  
  Thus we have  
  \begin{eqnarray*}
  f_C(q) & = & (1-q)(1-q\LL) \sum_{d=0}^\infty q^d [C^{[d]}] \\ 
  &  = & (1-q)(1-q\LL)  \sum_h [\{F\,|\,h_{F} = h\}] 
  \sum_{d=0}^{\infty} q^d [\proj^{h(d)-1}] \\
  & = & \sum_h [\{F\,|\,h_{F} = h\}] \cdot j_h(q)  .
  \end{eqnarray*}
  Collecting together the terms for $h$ and $h^\vee$ completes the proof. 
\end{proof}

For us it is more convenient to have the symmetry without the shifting. 

\begin{defn}
Let $C$ be a complete, reduced,  irreducible Gorenstein, complex curve of arithmetic genus $g$.  
Then we write
\[Z_C(q) := \sum_{n=0}^\infty  C^{[n]} q^{n+1-g}  .\]
Since $(\proj^1)^{[n]} = \proj^n$, we have
\[Z_{\proj^1}(q) = \frac{q}{(1-q)(1-q\LL)}\]
and 
define
classes $N^i_C \in \gring$ by the following formula: 
\[\sum_{n=0}^\infty  C^{[n]} q^{n+1-g} 
= \sum_{i=0}^\infty N^i_C \cdot Z_{\proj^1}^{i+1-g}.\]
\end{defn}

To be explicit, note that 
$Z:=\frac{q}{(1-q)(1-q\LL)}$ is a power series in $q(1+q\Z[\LL]\lbr q\rbr )$.  It therefore 
has a compositional inverse $q(Z) \in\Z[\LL]\lbr Z\rbr$ (given explicitly in footnote 6). Thus the $N^i_C$ are just defined by the change of variables
$$\sum_{i=0}^\infty N^i_C Z^{i+1-g}=\sum_{n=0}^\infty  C^{[n]} q(Z)^{n+1-g}.$$

\begin{cor} \label{cor:bps}
Let $C$ be a complete, reduced,  irreducible Gorenstein, complex curve of arithmetic genus $g$. 
Then $Z_C(q) = Z_C(1/q\LL)$, and $N^i_C = 0$ for $i > g$. 
\end{cor}
\begin{proof}
In terms of the $f_C$ of Proposition \ref{prop:ratfe}, we have 
$Z_C(q) =  q^{-g} Z_{\mathbb{P}^1}(q) f_C(q)$, so we may conclude 
the rationality and functional equation of $Z_C$ from that for $f_C$. 

Moreover $f_C(q)$ is a polynomial of degree $2g$.  On the other hand we
can expand it as 
$f_C(q) = q^g  \sum_{i=0}^\infty N^i_C \cdot Z_{\mathbb{P}^1}^{i-g}$
and then further into 
$$ \left(\sum_{i=0}^g N^i_C \cdot q^i (1-q)^{g-i} (1-q\LL)^{g-i}\right) + 
\left(q^g 
\sum_{i=1}^\infty 
 N^{g+i}_C \left( \frac{q}{(1-q)(1-q\LL)}\right)^{i}\right)\!.
$$
The first term is visibly a polynomial of degree at most 2g; for the second 
to be the same it must vanish. But the coefficient of $N^{g+i}_C$  in the second 
term is $q^{g+i}+ O(q^{g+i+1})$, so it follows that $N^{g+i}_C = 0$ for $i > 1$. 
\end{proof}

\begin{rem}  We always have $N^0_C = 1$ and $N^1_C = [C] + (g-1)(1+\LL)$. 
Given the vanishing, only $N^g_C$ contributes large
powers of $q$ to $Z_C(q)$; on the other hand when $n \gg 0$ the map 
$C^{[n]} \to \overline{J}^0(C)$
is a projective bundle.  Comparison of these terms reveals 
$N^g_C = \overline{J}^0(C)$. 
\end{rem}

\begin{example}
  Let $\P^1$, $A_1$, $A_2$ be rational curves
  that are smooth, have one node, and have one cusp respectively.  
  \begin{itemize}
    \item $N^0_{\P^1} = 1$. 
    \item $N^0_{A_1} = 1$ and 
      $N^1_{A_1} = \LL$.
    \item $N^0_{A_2}  = 1 $ and 
      $N^1_{A_2} = 1 + \LL$. 
  \end{itemize}
  In each case Corollary \ref{cor:bps} ensures
  the higher $N^i$ vanish.
\end{example}

As the Euler characteristic factors through $\gring$,  it makes
sense to write $\chi(N^i_C)$.  When $C$ is smooth of genus $g$, it follows
from Macdonald's calculation of the cohomology of symmetric products that
$\chi(N^i_C) = 0$ for $i > 0$.  However, 
$N^1_C = [C] + (g-1)(1+\LL)$ is never zero.  To avoid this incursion of
the global geometry of $C$, we remove the contribution of the normalization
$\widetilde{C}$ of $C$.  We define $\widetilde{Z}_C: = Z_C / Z_{\widetilde{C}}$, and
$\widetilde{N}^i_C$ by the formula:

\[
\widetilde{Z}_C 
= \sum_{i=0}^\infty \widetilde{N}^i_C
Z_{\proj^1}^{i-\delta}.
\]

For $p \in C$, let 
$C_p^{[n]}$ be the subvariety of $C^{[n]}$ parameterizing subschemes supported at $p$.  
Let $\delta(p)$ denote the $\delta$ invariant of the singularity at $p$, and $b(p)$ the number of 
analytic local
branches at $p$.  Note $C_p^{[n]}$ only depends on the analytic local structure of $C$ at $p$. 
A stratification argument gives the product expansion: 

\begin{equation} \label{eq:product}
 \widetilde{Z}_C =  \frac{\sum_{n=0}^\infty q^{n+1-g}[C^{[n]}]}{\sum_{n=0}^\infty q^{n+1-\widetilde g}[\widetilde C^{[n]}]} = 
\prod_{p \in C} \frac{\sum_{n=0}^\infty q^{n-\delta(p)} [C_p^{[n]}] }{(1-q)^{-b(p)}} .\end{equation}
The product is written over all points in $C$, but may as well be written over only singular 
points as the smooth points contribute $1$. 

We define $Z_{C, p}:= (1-q)^{b(p)} q^{-\delta(p)} \sum_{n=0}^\infty [C_p^{[n]}]$ in order that
the above can be written $\widetilde{Z}_C = \prod Z_{C, p}$.  We likewise define $N^i_{C, p}$ by
the expansion $Z_{C, p} :=  \sum_{i=0}^\infty N^i_{C, p} 
Z_{\proj^1}^{i-\delta(p)}$.

\begin{cor} \label{cor:vanishing} Let $C$ be a curve, let $p \in C$.  
Then $N^i_{C, p} = 0$ for $i > \delta(p)$, and 
$\widetilde{N}^i_C = 0$ for $i > \delta(C)$. 
\end{cor}
\begin{proof}
Let $C', p'$ be a rational curve with a unique singularity at $p'$ analytically isomorphic to the 
singularity of $C$ at $p$; such a curve exists by \cite[Prop.  2.1.1]{L}.  
Then
$$N^i_{C, p} =  N^i_{C', p'}  = \widetilde{N}^i_{C'}  = N^i_{C'}.$$  Since $\delta(p) = \delta(p') = g(C')$,
the vanishing  $N^i_{C,p}$ for $i > \delta(p)$ follows from Corollary \ref{cor:bps} applied to $C'$.  

Then we have 
$$\sum_{i=0}^\infty \widetilde{N}^i_C
Z_{\proj^1}^{i-\delta(C)} =   \widetilde{Z}_C = \prod_{p \in C} \sum_{i=0}^{\delta(p)} N^i_{C, p} 
Z_{\proj^1}^{i-\delta(p)} $$
giving the desired vanishing of $\widetilde{N}^i_C$ for $i > \delta(C)$. 
\end{proof}

\begin{conj} \label{conj:positive} 
  For all $h$, we have $\widetilde{N}^h_C \in \ZZ_{\ge 0}[\LL]$.  
\end{conj} 

\noindent {\bf Remark.}  Theorem \ref{thm:mult} realizes
$\chi(\widetilde{N}^i_C)$ as the multiplicity of the stratum 
of curves of cogenus $i$ inside the versal deformation of $C$, whence
it follows that $\chi(\widetilde{N}^i_C) > 0$ for $i \le \delta(C)$.  It may be hoped that
this conjecture indicates a refinement of this geometric structure,
i.e., that the coefficients of the $\widetilde{N}^i_C$ count something. 
\vspace{2mm}

For rational $C$, we have 
$\widetilde{N}^\delta_C = [\overline{J}(C)]$; by \cite{L} these Jacobians
are known to receive a bijective morphism from
an affine Springer fiber for $\mathfrak{gl}$.  It is known that
such affine Springer fibers in other types {\em are not necessarily} in
$\ZZ[\LL]$ \cite[Appendix]{KzL}, 
however according to Lusztig the status of the $\mathfrak{gl}$
affine Springer fibers is unknown.  
From the work of Piontkowski \cite{P} it follows that 
that for unibranch singularities with a single Puiseux pair, 
and for unibranch singularities
whose links are two-cablings of links of simple unibranch singularities,
one has at least $\widetilde{N}^\delta_C \in \ZZ[\LL]$.  
The stated positivity has been checked for unibranch 
singularities with a single Puiseux pair (e.g. $x^m = y^n$) for $m < 14$ and $n< 20$
using the explicit formula for $Z(C)$ given in \cite[Thm. 5]{ORSh}. 

\vspace{2mm}

From the fact that $\chi(Z_{\widetilde{C}}) = \chi(Z_{\proj^1})^{1-g}$, we see that $n^i_C :=
\chi(N^i_C) =  \chi(\widetilde{N}^i_C)$, and in particular that
the $n^i_C$ vanish for $i > \delta(C)$.  
This fact was used in \cite{PT2} as evidence that the $n^i_C$ 
were in fact the Gopakumar-Vafa invariants, and was 
exploited in \cite{KST} to count curves on surfaces. 
Here we note the Hirzebruch $\chi_{-y}$ genus has the same property, 
which suggests that it is a more sensible refinement than
working in the ring of varieties.

\begin{lem}
For $C$ a smooth curve of genus $g$, 
$\chi_{-y}(Z_C) = \chi_{-y}(Z_{\proj^1})^{1-g}$.  
\end{lem}
\begin{proof}
 The Hodge structures of symmetric
  products are known explicitly \cite{M}.
\end{proof}

\begin{cor}
$\chi_{-y}(N^i_C) = \chi_{-y}(\widetilde{N}^i_C)$.  In particular,
$\chi_{-y}(N^i_C) = 0$ for $i > \delta(C)$. 
\end{cor}

\begin{rem}
According to Conjecture \ref{conj:positive}, the $\widetilde{N}^i_C$ may
be recovered from $\chi_{-y}(N^i_C)$ by $y \mapsto \LL$.  Conversely, one
may separate Conjecture \ref{conj:positive} into two pieces, one asking that
$\widetilde{N}^i_C \in \Z[\LL]$ and the second asking that 
$\chi_{-y}(N^i_C)$ have positive coefficients. 
\end{rem}

\subsection{Curves with simple singularities} \label{subsec:simple}

\begin{figure}
\framebox{\includegraphics[width=14cm]{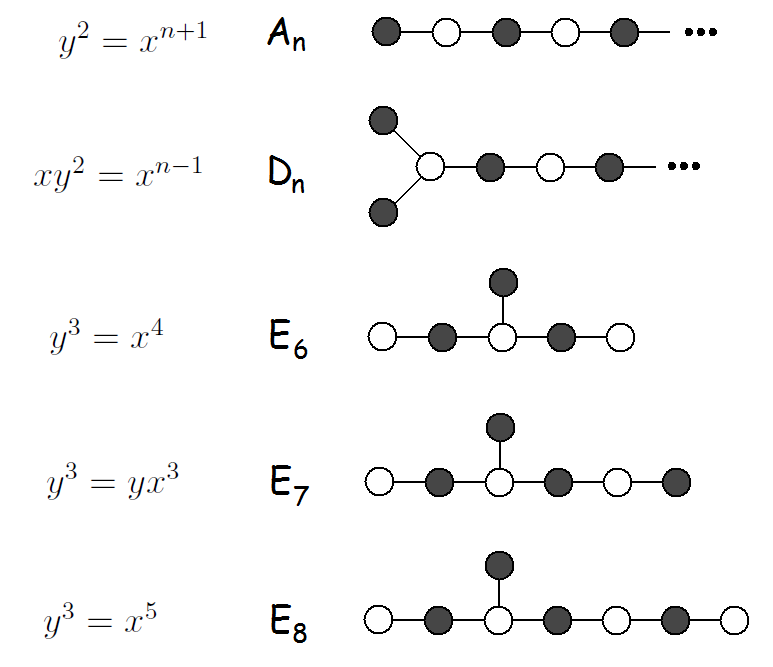}}
\caption{ \label{fig:ADE} The ADE singularities and associated colored diagrams.
The subscript gives the total number of vertices of the diagram, and
the Milnor number of the singularity.  
The coloring is characterized by requiring the colors to alternate, 
and requiring that the total number of black vertices is the delta invariant. 
Except for $A_1$, the number of filled vertices of valence one is the number of analytic
local branches.}
\end{figure}

Here we express the refined BPS invariants of curves with simple singularities
in terms the associated Dynkin diagram. 
In Figure \ref{fig:ADE}, we recall the ADE classification of simple
singularities. 

\begin{lem} \label{lem:ADEinfinity}
We define $A_\infty$ to be the germ at the origin of the 
curve cut out by $y^2 = 0$.  Similarly 
we define $D_\infty$ by $xy^2 = 0$ and $E_\infty$ by $y^3 = 0$.  
Then, for $X = A,D,E$, we have an equality
$X_\mu^{[i]} = X_\infty^{[i]}$ as subsets of $(\CC^2)^{[i]}$ 
for any $i$ up to the delta invariant of $X_\mu$. 
\end{lem}
\begin{proof}
Any subscheme of length at most $\delta$ supported at the 
origin is annihilated by $(x,y)^\delta$.  In each case, the RHS of
the equation of $X_\mu$ already belongs to this ideal.
\end{proof}

\begin{prop} \label{prop:ADEinfinity}
\begin{eqnarray*}
  \sum [A_\infty^{[n]}] q^n & = & \frac{1}{(1-q)(1-q^2 \LL)}, \\
  \sum [D_\infty^{[n]}] q^n & = & \frac{1-q+q^3 \LL^2}{(1-q)^2(1-q^2 \LL)}, \\
  \sum [E_\infty^{[n]}] q^n & = & \frac{1}{(1-q)(1-q^2 \LL)(1-q^3 \LL^2)}. \\  
\end{eqnarray*}
\end{prop}
\begin{proof}
%
  We fix the  monomial order $1 < x < x^2 < \cdots < xy < x^2 y < x^3 y  < \cdots$.  Recall that in the setting of power series
  rings, the theory of Gr\"obner bases (or `standard bases') is developed with respect to the {\em lowest} rather than highest order
  terms of series.  Thus let $lt(f)$ of 
  $f \in \CC\lbr x, y\rbr$ be the lowest degree monomial appearing.  A generating set 
  $I = (i_1, i_2, \ldots, i_k)$ is a Gr\"obner basis when $lt(I) = (lt(i_1), \ldots, lt(i_k))$.  It's reduced when 
  the $lt(i_t)$ are a minimal generating set for $lt(I)$, i.e. when $k = \dim_\CC lt(I) / (x,y) lt(I)$, and moreover 
  $$i_t = lt(i_t) + \mbox{monomials not in $lt(I)$}. $$
  As is well known, every ideal admits a unique reduced Gr\"obner basis. 
  

  Let's first consider the case $A_\infty$.  Ideals of $\CC\lbr x, y\rbr / y^2$ are the same as ideals
  of $\CC\lbr x, y\rbr$ containing $y^2$.  The lt ideal of a finite colength such ideal must be of the form
   $(y^2, x^a)$ or $(y^2, yx^{b}, x^a)$ for some $b < a$.  In the first case, a reduced Gr\"obner basis
  of the original ideal takes the form $$\left(y^2, x^a + y \sum_{i=0}^{a-1} c_i x^i \right)\! ,$$
  whereas in the second it  takes the form
  $$\left(y^2, yx^{b}, x^a + y \sum_{i=0}^{b-1} c_i x^i \right)\! .$$
  These sum as: 
  $$\sum [A_\infty^{[n]}] q^n = \sum_{a=0}^\infty \sum_{b = 0}^a q^{a+b} \LL^b = 
   \sum_{a=0}^\infty q^a \frac{1-(q\LL)^{a+1}}{1-q\LL} = \frac{1}{1-q\LL} 
   \left(\frac{1}{1-q} - \frac{q\LL}{1-q^2 \LL} \right)\!,$$
  which simplifies to the stated expression.  
  The $E_\infty$ case can be treated similarly, and in any event the statements for $A_\infty, E_\infty$ are special cases of \cite[Prop. 6]{ORSh}. 
  
  We turn to $D_\infty$.  Let $I \subset \CC\lbr  x , y\rbr$ be an ideal containing $xy^2$ and of finite colength.  We have already
  counted all the ideals containing $y^2$; these contribute $1/(1-q)(1-q^2 \LL)$. 
  We must account for the remaining ideals, and show
  they contribute $q^3 \LL^2 / (1-q)^2(1-q^2 \LL)$. 
  So let $I \subset \CC\lbr x, y \rbr$ be an ideal containing $xy^2$ but not $y^2$.  Then $lt(I)$ contains $xy^2$; it may
  be written in the form $(y^a, xy^2, x^b y, x^c)$ for some $b \le c$.  A corresponding Gr\"obner basis is
  $$\left(y^a, xy^2, x^b y + \sum_{i=2}^{a-1} c_i y^i, x^c + y \sum_{j=1}^{b-1} d_j x^j + \sum_{k=1}^{a-1} e_k y^k \right)\! .$$
  
  We start the sum over $j$ from $1$ since
  the monomial $y$ is already accounted for in the sum over $k$. 
  Note that $a \ge 3$  and $b, c \ge 1$ because we assumed $y^2 \notin I$.   There are 
  further constraints on which such things may be Gr\"obner bases.  
  For instance, multiplying the third generator by $y$ and subtracting a multiple of $xy^2$ gives
  $\sum c_i y^{i+1}$; by our assumption on $lt(I)$ lowest term of this must be $c_{a-1} y^a$, i.e., the other $c_i$ must vanish. 
  So the basis takes the form
  $$\left(y^a, xy^2, x^b y +   c' y^{a-1}, x^c + y \sum_{j=1}^{b-1} d_j x^j + \sum_{k=1}^{a-1} e_k y^k \right)\!.$$
  Multiplying the fourth generator by $y^2$ and arguing similarly, we see in fact $e_k$ vanishes for $k < a-2$, so
  the basis takes the form
  $$\left(y^a, xy^2, x^b y +   c' y^{a-1}, x^c + y \sum_{j=1}^{b-1} d_j x^j + e y^{a-2} + e' y^{a-1} \right)\!.$$
  Multiplying the fourth element by $y$ and subtracting from the third element times
  $x^{c-b}$ gives $ey^{a-1} + e'y^a -  c' x^{c-b} y^{a-1} = 
  y^{a-1}(e + e'y - c' x^{c-b})$.  If $e \ne 0$, then the term in parenthesis is invertible and 
  hence $y^{a-1}$ is in the ideal, a contradiction.  Finally we are reduced to the form
  $$\left(y^a, xy^2, x^b y +   c' y^{a-1}, x^c + y \sum_{j=1}^{b-1} d_j x^j + e' y^{a-1} \right)\!.$$

  We leave it to the reader to check that there are no other constraints, and moreover that 
  every ideal admits a unique basis of this form (even though it is not literally a reduced
  Gr\"obner basis, e.g. when $b = c$).  We count them: 
  
  $$\sum_{a=3}^\infty \sum_{c=1}^\infty \sum_{b=1}^c \LL^{b-1+2} q^{a+b+c - 2} =
  q^3 \LL^2 \sum_{a' = 0}^\infty \sum_{c'=0}^\infty \sum_{b'=0}^{c'} \LL^{b'} q^{a'+b'+c'} 
  = \frac{q^3 \LL^2}{(1-q)^2 (1-q^2 \LL)}.$$
\end{proof}

Let $X$ be a simple singularity type.  If $C$ is a curve and $p$ a point at which $C$ has a singularity
analytically of type $X$, we write $X^{[n]} := C_p^{[n]}$, and similarly
$\delta(X)$ and $b(X)$ for what we have before written as $\delta(p)$ and $b(p)$; 
and similarly $Z_X := Z_{C,p}$ and $N_X^i := N_{C, p}^i$.  Recall these are related by:

$$\sum_{i=0}^{\delta(X)} N_X^i \cdot Z_{\P^1}^{i - \delta} = 
Z_X = (1-q)^{b(X)} q^{-\delta(X)} \sum_{n=0}^\infty q^n X^{[n]}.$$

\begin{thm} \label{thm:simple}
  Let $X$ be a simple singularity.  
  Color the dots of the
  associated Dynkin diagram as in Figure \ref{fig:ADE}.   
  Let $n^{w,b}_X$ be the number
  of ways to choose $w$ white dots and $b$ black ones such that no
  two dots 
  are adjacent.  Then 
  \[N^h_X = \sum_{w+b=h} \!n^{w,b}_X\, \LL^b. \] 
\end{thm}
\begin{proof}
  We temporarily write 
  \[M^h_X := \sum_{w+b=h} \!n^{w,b}_X\, \LL^b. \]
  and $Y_X :=  \sum_{i=0}^{\delta(X)} M_X^i \cdot Z_{\P^1}^{i - \delta}$.
  We will show $Y_X = Z_X$.  For small cases, say 
  $X = A_1, A_2, A_3 = D_3, D_4, E_6, E_7, E_8$, we have explicit formulas
  on both sides so this may be verified by hand or by computer.  
  
  It remains
  to treat in general $A_n, D_n$.  The argument  will reveal that
  in fact, Lemma \ref{lem:ADEinfinity} in a sense determined the values of the series in 
  Proposition \ref{prop:ADEinfinity}; in particular, we will not use that Proposition again. 
  (We used it above to check the base cases; of course these could have been done 
  without appealing to its full strength.) 
    

  We begin with the case of $A_n$.  
  By considering what happens according as the 
  right end dot of the Dynkin diagram -- the one not seen in the pictures above,
  which is white if $n = 2k$ and black if $n = 2k-1$ -- is chosen or not, we have 
  $$n_{A_{2k}}^{w, b} = n_{A_{2k-1}}^{w,b} + n_{A_{2k-2}}^{w-1,b} \qquad \qquad \qquad
  n_{A_{2k-1}}^{w, b} = n_{A_{2k-2}}^{w,b} + n_{A_{2k-3}}^{w,b-1}.$$
  Summing these, we have 
  $$M_{A_{2k}}^h = M_{A_{2k-1}}^h + M_{A_{2k-2}}^{h-1} \qquad \qquad \qquad
  M_{A_{2k-1}}^h = M_{A_{2k-2}}^h + \LL M_{A_{2k-3}}^{h-1}.$$
  Noting that $\delta(A_{2k}) = \delta(A_{2k-1}) = k$ and summing,
  $$Y_{A_{2k}} = Y_{A_{2k-1}} + Y_{A_{2k-2}} \qquad \qquad \qquad
  Y_{A_{2k-1}} = Y_{\P^1}^{-1} Y_{A_{2k-2}} + \LL Y_{A_{2k-3}}.$$
  It remains to show the same for the $Z$'s.   
  Since both are $Y_X$ and $Z_X$ are symmetric Laurent polynomials in $q$, 
  it suffices to check the recurrence holds mod $q$.  Writing $\xi_A = \sum q^n A_\infty^{[n]}$
  we had from Lemma \ref{lem:ADEinfinity} that 
  $Z_{A_n} = (1-q)^{b(A_n)} q^{-\delta(A_n)} \xi_A \pmod q$.  We are left to show
  \begin{eqnarray*} q^{-k} (1-q) \xi_A & =  & q^{-k}(1-q)^2 \xi_A + q^{-k+1} (1-q) \xi_A \pmod q, \\
  q^{-k}(1-q)^2 \xi_A & = & q^{-k} (1-q)^2(1-q\LL) \xi_A + \LL q^{-k+1} (1-q)^2 \xi_A \pmod q.
  \end{eqnarray*}
  But these are both formal identities, independent of the value of $\xi_A$.  
  
  An identical argument suffices to treat the $D_n$.  We remark in passing that
  the recursion above corresponds \cite{OSh, ORSh, Mau} to the skein relation
  on the HOMFLY invariants of the links of the singularities. 
  
\end{proof}

\begin{example}  For $A_9$, the nonvanishing invariants are
  \begin{eqnarray*}
    \widetilde{N}^4 & = & 1 + \LL + \LL^2 + \LL^3 + \LL^4 \\
    \widetilde{N}^3 & = & 4 + 6 \LL + 6 \LL^2 + 4 \LL^3 \\
    \widetilde{N}^2 & = & 6 + 9 \LL + 6 \LL^2 \\
    \widetilde{N}^1 & = & 4 + 4 \LL \\
    \widetilde{N}^0 & = & 1
  \end{eqnarray*}
\end{example}

The $\widetilde{N}^i$ are not generally symmetric. 

\begin{example} For $E_6$, the nonvanishing invariants are
  \begin{eqnarray*}
    \widetilde{N}^3 & = & 1 + \LL + 2 \LL^2 + \LL^3 \\
    \widetilde{N}^2 & = & 3 + 4 \LL + 3 \LL^2 \\
    \widetilde{N}^1 & = & 3 + 3 \LL \\
    \widetilde{N}^0 & = & 1 \\
  \end{eqnarray*}
\end{example}

\begin{example} For $E_8$, the nonvanishing invariants are
  \begin{eqnarray*}
    \widetilde{N}^4 & = & 1 + \LL + 2\LL^2 + 2\LL^3 + \LL^4 \\
    \widetilde{N}^3 & = & 4 + 6 \LL + 7 \LL^2 + 4 \LL^3 \\
    \widetilde{N}^2 & = & 6 + 9 \LL + 6 \LL^2 \\
    \widetilde{N}^1 & = & 4 + 4 \LL \\
    \widetilde{N}^0 & = & 1
  \end{eqnarray*}
\end{example}

For the simple singularities, we have the following remarkable statement, 
which may be proven
by comparing Theorem \ref{thm:simple} to the description of the versal
deformation of a simple singularity as the quotient of the hyperplane arrangement
of the same name by the Weyl group \cite{AGV}. 

\begin{thm}\label{duco}
  Let $X$ be a simple singularity.  Then there exists some curve $C$ containing
  as its unique singularity a real form of $X$, and a real disc $D^j$ in the real
  locus of the versal deformation of $X$ such that 
  \[\chi_{-y}(N^j_C) = \sum_i y^i \cdot \# (D^j \cap B^{j, i}_+),\]
  where $B^{j,i}_+$ is the locus of real nodal curves with $j$ total nodes of which $i$ 
  are hyperbolic.
\end{thm}

This result was known to van Straten when $j = \delta(C)$; he 
had conjectured in this case that it holds for all singularities (see \cite{vS}[Conj.~4.7]) on the evidence
of its validity for the simple singularities and its validity at $y = 1$ for all singularities \cite{FGvS}.  
Theorem \ref{thm:mult} asserts the validity of the above statement at $y=1$ for all singularities, so one
might analogously conjecture that the statement of Theorem \ref{duco} holds always. 
Note that for unibranch singularities, there is a unique topological type of real form. 

\subsection{Real invariants}

We recall that the Betti numbers may depend on the choice of coefficients for cohomology: for instance, we have 
$$\sum q^n \dim \mathrm{H}^k (\P^2(\R), \Q) = 1 \ne 1 + q + q^2 = \sum q^n \dim \mathrm{H}^k (\P^2(\R), \Z/2\Z).$$
Note that nonetheless we have $\chi(\P^2(\R), \Q) =1 = \chi(\P^2(\R), \Z/2\Z)$.  In fact, as is well known to follow
from the universal coefficient theorem, the Euler
characteristic does not depend on the choice of coefficients.   The same holds for compactly supported cohomology,
for instance because this can be described as the relative cohomology of an appropriate pair.  We therefore suppress
the coefficients from the notation for Euler characteristics.

In particular, since any (not necessarily oriented or closed) manifold $M$ enjoys Poincar\'e duality between 
$\mathrm{H}^*(M, \Z/2\Z)$ and $\mathrm{H}^*_c(M, \Z/2\Z)$, we have $\chi_c(M) = (-1)^{\mathrm{dim}\, M} \chi(M)$.  

As pointed out by Macdonald \cite{M}, for any reasonable space $X$, we have 
$$\mathrm{H}^*(X^n/\mathfrak{S}_n, \Q) = \mathrm{H}^*(X^n,\Q)^{\mathfrak{S}_n} = 
\left(\mathrm{H}^*(X, \Q)^{\otimes n}\right)^{\mathfrak{S}_n}$$
and it follows formally that
$$\sum_{n=0}^\infty q^n \chi(X^n/\mathfrak{S}_n) = \left(\frac{1}{1-q}\right)^{\chi(X)}.$$ 
We note this is false in general for the compactly supported Euler characteristic, for instance,  
one can show
$\sum_{n=0}^\infty q^n \chi_c(\R^n/\mathfrak{S}_n) = 1$. 
Nonetheless, for real surfaces we have:

\begin{lem} \label{lem:2man}
Let $\Sigma$ be a manifold without boundary of (real) dimension 2; it may be non-orientable and non-compact. 
Then 
$$\sum_{n=0}^\infty q^n \chi_c(\Sigma^n / \mathfrak{S}_n )  =  \left(\frac{1}{1-q}\right)^{\chi_c(\Sigma)}.$$
\end{lem}
\begin{proof}
Both $\Sigma^n/\mathfrak{S}_n$ and $\Sigma$ are real even dimensional manifolds, so have $\chi = \chi_c$. 
\end{proof}


Now let $C$ be an algebraic curve defined over $\mathbb{R}$.  We consider

\[Z^\R_C:= \sum_{n=0}^\infty q^{n+1-g} \chi_c(C^{[n]}(\R)).\]

Note there are {\em two} different projective lines over $\R$: one which has a real point (and hence whose real locus is a circle),
say $X^2 - Y^2 + Z^2 = 0$,  and one which has no real points $X^2 + Y^2 + Z^2 = 0$.  We call them $\P^1_-$ and $\P^1_+$.  Their symmetric
powers likewise carry different complex conjugations: 
$$(\P^1_-)^{[n]} (\R) = \R \P^n \,\,\,\,\,\,\,\,\,\,\,\,\,\,\,\,\,\,\,\,\,\,\,\,\,\,\,\,\,\, 
(\P^1_+)^{[2k+1]}(\R) = \emptyset \,\,\,\,\,\,\,\,\,\,\,\,\,\,\,\,\,\,\,\,\,\,\,\,\,\,\,\,\,\, 
(\P^1_+)^{[2k]}(\R) = \R \P^{2k}. $$ 
Nonetheless, since $\chi_c(\R \P^{2k+1}) = 0$, we have 
$$Z_{\P^1_+}^\R = q/(1-q^2) = Z_{\P^1_-}^\R.$$

In fact, for a smooth curve $C$, the series $Z_C^\R$ does not depend on the real structure of $C$. 
We can determine it explicitly:

\begin{lem} \label{lem:realbpssmooth}
Let $C$ be a curve defined over $\R$ such that $C \otimes_\R \C$ is a smooth curve of genus $g$.  Then
$Z_C^\R = (Z^\R_{\proj^1})^{1-g}$.
\end{lem}
\begin{proof}
   As $C$ is smooth, the Hilbert schemes agree with the symmetric products.  Let $\sigma$ denote 
   the complex conjugation.  We may stratify the locus $(\mathrm{Sym}^n(C)) (\R)$ according
   to the number of pairs of complex conjugate points.  We parameterize the real points by the symmetric powers
   of the real locus, and $n$ pairs of complex conjugate points in
   $C$ by $\mathrm{Sym}^n \left( \frac{C(\C) \setminus C(\R)}{\sigma} \right)$.  Thus compatibility of compactly
   supported Euler characteristic with cut-and-paste gives
  \[\sum_{n=0}^\infty q^n \chi_c \left( \mathrm{Sym}^n(C) (\R) \right) = 
  \left(\sum_{n=0}^\infty q^n \chi_c \left( \mathrm{Sym}^n(C(\R)) \right) \right)
   \left(\sum_{n=0}^\infty
    q^{2n} \chi_c \left( \mathrm{Sym}^n \left( \frac{C(\C) \setminus C(\R)}{\sigma}  \right) \right) \right)\!.\]
  The loci $ \mathrm{Sym}^n(C(\R)) $ are products of symmetric products of the connected components
  of $C(\R)$, i.e., of circles; it is easy to see (for instance by 
  using a circle action) that these have Euler characteristic and compactly supported Euler characteristic zero. 
  On the other hand, the second term was computed above in  Lemma \ref{lem:2man}.   We conclude
  \[
  \sum_{n=0}^\infty q^n \chi_c \left( \mathrm{Sym}^n(C) (\R) \right) = 
   \sum_{n=0}^\infty
    q^{2n} \chi_c \left( \mathrm{Sym}^n \left( \frac{C(\C) \setminus C(\R)}{\sigma}  \right) \right) 
     =  \left( \frac{1}{1-q^2} \right)^{\chi_c\left(\frac{C(\C) \setminus C(\R)}{\sigma} \right)}.\] 
  Finally, since a circle has vanishing Euler characteristic, and since the Euler characteristic
  of a $n:1$ \'etale cover is just $n$ times the Euler characteristic of the base, we have 
  $\chi_c\left(\frac{C(\C) \setminus C(\R)}{\sigma} \right) = \chi(C(\CC)) / 2 = 1-g$. 
\end{proof}

This motivates the definition of integers
$n^{i,\R}_C$ by the following formula:

\begin{equation}
\sum_{n=0}^\infty \chi_\R(C^{[n]})  \,\, q^{n-g+1} = 
\sum_{i=0}^\infty n^{i,\R}_C \left(\frac{q}{1-q^2} \right) ^{i-g+1} .
\end{equation}

In this notation, Lemma \ref{lem:realbpssmooth} asserts that, for $C$ a smooth
curve defined over $\R$, we have $n^{0,\R}_C = 1$ and $n^{i,\R}_C = 0$ for $i > 0$.

For a point $p$ on a curve $C$, write $(C,p)$ for the germ at $p$, and 
$(C,p)^{[n]}$ for the locus in the Hilbert scheme of points on $C$ of
subschemes set-theoretically supported at $p$.  For $p \in C(\R)$,
we write $rb(p)$ and
$cb(p)$ for respectively the number of real and complex points above $p$
in the normalization of $C$. 

We have by stratification
\begin{eqnarray*}
  \frac{\sum q^{n+1-g} \chi_c (C^{[n]}(\R))  }{\sum q^{n+1-\widetilde{g}} 
    \chi_c (\widetilde{C}^{[n]}(\R))} & = &
\prod_{p \in C^{\mathrm{sing}}(\R)}q^{-\delta(p)}\frac{\sum q^n \chi_c  ( (C,p)^{[n]}(\R) )}
{(1-q^2)^{-cb(p)/2} (1-q)^{-rb(p)}} \\ & \times &
\prod_{p \in C^{\mathrm{sing}}(\CC)/\sigma}q^{-2\delta(p)} 
\frac{\sum q^{2n} \chi_c ( (C,p)^{[n]})} 
{(1-q^2)^{cb(p)}}.
\end{eqnarray*}

We now compute $n^{\delta,\R}_C$ for a curve $C$ 
with $\delta$ nodes and no other singularities. 
By the product formula above, it suffices to do this 
in three special cases.

\begin{lem}
  The only nonvanishing invariants
  of a curve $c_{\pm}$ of arithmetic genus $1$ with a single node
  analytically of the form $\R\lbr x,y\rbr /(x^2 \pm y^2)$ 
  are 
  $n^{0,\R} = 1$ and $n^{1,\R} = \pm 1$. 
\end{lem}

(Caution: $c_-$ is the one with a node that looks like $+$.) 

\begin{lem}
  The only nonvanishing invariants of 
  a curve of arithmetic genus two with a pair of complex
  conjugate nodes are 
  $n^{0,\R}  = 1$ and $
  n^{1,\R}  = 0$ and $n^{2,\R} = 1$. 
\end{lem}

We conclude:

\begin{prop} 
For a nodal curve $C$ with $\delta = \delta_+ + \delta_- + 2 \delta_0$ nodes, 
where $\delta_{\pm}$ are of the form $\R\lbr x,y\rbr /(x^2 \pm y^2)$
and the $2 \delta_0$ are complex conjugates, 
\[n^{\delta,\R}_C = (-1)^{\delta_-} = (-1)^{\delta_+ - \delta}.\]  
Moreover, $n^{i,\R}_C = 0$ for $i > \delta$. 
\end{prop}

\begin{thm} \label{thm:realmult}
  Let $C$ be a real reduced plane curve, and let $\cC \to B$ be a 
  versal deformation of its singularities.  Let $B^{\delta, \delta_-}_+ \subset B(\R)$
  be the locus of nodal curves with $\delta$ nodes of which $\delta_-$ 
  are of the form $\R\lbr x,y\rbr /(x^2 - y^2)$.  Let $D^j$ be a general disc of dimension $j$, preserved by complex conjugation,
  passing near $[C] \in B$.
  Then \[n^{j,\R}_C = \sum_k (-1)^k D(\R)^j \cap B_+^{j, k}.\]
  In particular, $n^{j,\R}_C = 0$ for $j > \delta(C)$. 
\end{thm}
\begin{proof}
In the following we write $\delta:=\delta(C)$.
   View $n^{j,\R}$ as a constructible function on $B(\R)$ taking
   $b \mapsto n^{j,\R}_{\cC_b}$.

   Let $D^i_0 \subset B$ be a general (complex but preserved
   by conjugation) disc of dimension $i$ containing $[C]$.  Then according to \cite{FGvS, Sh}, the
   first $\cC^{[\le i]}_{D^i_0}$ are all smooth, and if $i \ge \delta$ then all $\cC^{[n]}_{D^i_0}$ are smooth. 
   The real locus of a smooth variety is smooth, so 
   the same holds upon passing to real points. Taking $D^{\delta+1}$ a disc containing $D^\delta_0$ and
   $D^\delta_1$ a sufficiently nearby slice, we have by smoothness that the spaces
   $\cC^{[n]}_{D^\delta_0}(\R)$ and $\cC^{[n]}_{D^\delta_1}(\R)$ are diffeomorphic and therefore have the same 
   compactly supported Euler characteristics.  Note that the first $j$ Hilbert schemes suffice to 
   determine $n^{j,\R}$.  By additivity it follows that
   \begin{equation}
   \label{eq:intequal} 
   \int_{D^i_0(\R)} n^{j,\R} d\chi = \int_{D^i_1(\R)} n^{j,\R} d \chi 
   \,\,\,\,\,\,\,\,\,\,\,\,\mbox{for any $j \le i$, and for any $j$ at all if $i \ge \delta$.}
   \end{equation}

   We first show the vanishing of the $n^{j,\R}$ for $j > \delta(C)$.  Note  we have already
   established it for smooth and nodal curves.  We induct on $\delta(C)$.  Take $i = \delta$ in 
   Equality \ref{eq:intequal}.
   By \cite{DH,T}, the locus of curves of cogenus at least $\delta$ is of codimension $\delta$ and is the closure of the 
   locus of $\delta$-nodal curves.  Thus by genericity its 
   only intersection with $D^\delta_0$ is at the central point $[C]$, and its only intersection with $D^\delta_1$
   is in finitely many nodal curves.  Thus by induction 
   and our explicit verification in the case of nodal curves, 
   the integral on the right vanishes for $j > \delta$, 
   hence so does the integral on the left, which again by 
   induction is equal to $n^{j,\R}_C$. 
   
   Now we consider the remaining $n^j$.  Take $i = j$ in the above equality.  Then by the same reasoning the only
   contribution to the integral on the left is $n^{i,\R}_C$, and the only contribution to the integral on the right
   is the $j$-nodal curves.  By our previous calculation, these contribute as required.
  \end{proof}

\begin{rem}
	If $C$ has a real line bundle of degree $1$ then the proof
	of Proposition \ref{prop:ratfe} and Corollary \ref{cor:bps} carries through in $\mathrm{K}_0(\mathrm{var}/\R)$; 
	we may use this instead to conclude that $n^{i,\R}_C = 0$ for $i > g(C)$.  However
	we have not found an analogous argument if $C$ has no such bundle. 
\end{rem}

\begin{rem}
        In \cite{Sh}, the proof of Theorem \ref{thm:mult} used as an input the vanishing of the 
        $n^i_C$ for $i > \delta$; the argument given above shows this was unnecessary.
\end{rem}

\begin{rem}
	Thus a necessary condition for the statement of Theorem \ref{duco} to hold for a singularity
	$X$
	is the existence a curve $C$ containing a real form of $X$ with $\chi_{1}(Z_C) = Z^\R_C$.  
\end{rem}

\begin{cor} \label{cor:realbps}
  Let $\cC \to B$ be a family of reduced real plane curves 
  in which all curves have cogenus $\le \delta$ 
  and there are finitely many curves of cogenus $\delta$, all nodal.  Then
  defining $n^{i,\R}_{\cC/B}$ by   
  \[
  \sum_{n=0}^\infty \chi_c(\cC_B^{[n]}(\R))  \,\, q^{n-g+1} = 
\sum_{i=0}^\infty n^{i,\R}_{\cC/B} \left(\frac{q}{1-q^2} \right) ^{i-g+1} .
\]
  the number of $\delta$-nodal curves counted with signs 
  $(-1)^{\delta_-} = (-1)^{\delta_+ - \delta}$ is $n^{\delta,\R}_{\cC/B}$.
\end{cor}

\section{Refined invariants of linear systems}
\label{sec:families}

Let $\pi:\cC \to B$ be a family of plane curves, and $\pi^{[n]}:\cC^{[n]}\to B$
the relative Hilbert schemes.  Denote by $\mhm(B)$ the category of mixed Hodge modules \cite{Sai}
over $B$.  We define:
\[
\mathcal{Z}_{\cC/B} = \sum_{n=0}^\infty \pi^{[n]}_! \Q_{\cC^{[n]}} q^{n+1-g} \in K_0(\mhm(B))\lbr q\rbr .
\]
We define invariants
$\mathcal{N}^i_{\cC/B} \in K_0(\mhm(B))$ by the formula
\begin{equation} \label{eq:nyMHM}
\sum_{n=0}^\infty \pi^{[n]}_! \Q_{\cC^{[n]}} q^{n+1-g} = 
\sum_{i=0}^\infty  \mathcal{N}^i_{\cC/B} \times \mathcal{Z}_{\P^1/\mathrm{pt}}^{i+1-g} .
\end{equation}

\begin{prop} \label{prop:mhmbps}
  If $\pi:\cC \to B$ is a family of integral plane curves and moreover 
  $\cC^{[n]}$ is smooth for all $n$, then $\mathcal{N}^i_{\cC/B}=0$ for $i > g$.
\end{prop}
\begin{proof}
  In fact, according to \cite{MY,MS} we know much more.  Writing $\widetilde{\pi}$ for 
  the restriction of the map to the locus on the base 
  where it is smooth, and  $(1)$ for the Tate twist (it decreases weights by $2$, so $[2](1)$ preserves weights)
  we have:   
  \[(1-q)(1-q \Q_B[2](1))  \sum_{n=0}^\infty q^{n} \pi^{[n]}_! \Q_{\cC^{[n]}}[\dim B + n] = 
    \sum_{n=0}^{2g} q^n \mathrm{IC}(B, \Lambda^i R^1 \widetilde{\pi}_* \Q_{\cC}[\dim B + g])[-i].
  \]
  The RHS is a polynomial of degree $2g$ which enjoys the symmetry
  $\Lambda^i R^1 \cong \Lambda^{2g-i} R^1$ by hard Lefschetz; we  
   deduce the vanishing of the $\mathcal{N}^i_{\cC/B}$ as in Corollary \ref{cor:bps}.   
\end{proof}

\begin{rem}
  One may write the same definition in $\mathrm{K}_0(\mathrm{var}/B)$, but then 
  we do not know whether $\mathcal{N}^i$ vanishes for $i > g$.  The problem
  arises already for families of smooth curves; the question here is
  whether a family of Jacobians is equivalent to its torsors in the Grothendieck 
  group of varieties.  
\end{rem}

We may take pointwise Euler characteristic to define an
integer valued constructible function
$n^i:=\chi(\mathcal{N}^i)$, or global Euler characteristic 
to define $n^i_{\cC/B} \in \ZZ$.  
(By proper base-change, these are the same 
$n^i$ as in the introduction.) Since $n^i$ is supported on the locus
of curves of cogenus $i$, certainly $n^k_{\cC/B}$ vanishes if $k$ is greater
than the maximum cogenus of any curve in the family.  These are the same
as the constructible functions of the introduction, that is, they satisfy (and 
could take as their definition) the formula 
\[\sum_{n=0}^\infty \chi(\cC^{[n]}) q^{n+1-g} = 
\sum_{i=0}^\infty  n^i_{\cC/B} \left( \frac{q}{(1-q)^2}\right)^{i+1-g}.\] 

\begin{defn}\label{def:change}
We define
$N^i_{\cC/B}(y):=\chi_{-y}((B \to pt)_! \mathcal{N}^i)$.  Equivalently
we can apply $\chi_{-y} \circ (B \to pt)_!$ to both
sides of Equation \ref{eq:nyMHM} in order to directly define $N^i_{\cC/B}(y)$
by the formula \[\sum_{n=0}^\infty \chi_{-y}(\cC^{[n]}) q^{n+1-g} = 
\sum_{i=0}^\infty  N^i_{\cC/B}(y) \left( \frac{q}{(1-q)(1-qy)}\right)^{i+1-g}.\] 
\end{defn}

From Proposition \ref{prop:mhmbps}, we see that $N^i_{\cC/B}(y)$  vanishes  for 
$i$ greater than the maximum genus of any curve in a family whose relative Hilbert schemes are nonsingular.  
But for both the $\chi_{-y}$ invariants of a single curve, and the $\chi$ invariants of families,
we had vanishing beyond the maximum {\em cogenus}. 
Thus we may at least plausibly ask whether this holds for $N^i_{\cC/B}$.  In fact it need not:
Fantechi and Pandharipande
observed that this vanishing can fail already for $B$ a curve of positive genus 
and $\cC \to B$ a family of smooth curves.  However, empirically, 
the situation appears to be better for linear systems of curves in surfaces.
We have the following conjecture:

\begin{conj}\label{cochyfam}  Let $L$ be a line bundle on a surface $S$, 
$\cC \to |L|$ the tautological family of curves, and
$\CP^\delta \subset |L|$ be a linear subsystem. 
Assume the relative
Hilbert schemes $\cC^{[n]}_{\CP^\delta}$ are nonsingular for all $n \ge 0$. 
Then $N^i_{\cC/\CP^\delta}(y) = 0$ for $i > \delta$. 
\end{conj}

From the smoothness
criterion in \cite{FGvS} one may deduce that the maximum cogenus of any curve in any
such $\proj^\delta$ is $\delta$.  The assumption holds in the following situations: 

\begin{thm} \label{thm:whensmooth}
Let $L$ be a line bundle on a surface $S$,
$\cC \to |L|$ the tautological family of curves, and $\proj^\delta \subset |L|$ a general linear subsystem.  
Then all relative Hilbert schemes $\cC^{[n]}_{\CP^\delta}$ are nonsingular in the following situations: 
\begin{itemize}
	\item $S$ is arbitrary and $L$ is $\delta$-very ample \cite{KST, FGvS, Sh}.\footnote{
	What is actually observed in \cite{KST} is the (obvious) fact that the first $\delta$ relative Hilbert
	schemes are nonsingular under this hypothesis, and the  less obvious fact that the assumption 
	implies that no non-reduced curves or curves of cogenus $> \delta$ occur in the $\proj^\delta$.  
	Consequently the smoothness criteria in \cite{FGvS, Sh} may be used to establish smoothness
	of the remaining relative Hilbert schemes. 
	}
	\item $S$ is a K3 or abelian surface and $L$ is irreducible \cite{Muk}.
	\item $S$ is a rational surface and $\proj^\delta$ contains no non-reduced curves, and
	no curves with components which intersect $K_S$ non-negatively \cite{KS}.  In particular,
	for a general $\proj^{\le 2d-2} \subset |\mathcal{O}_{\proj^2}(d)|$. 
\end{itemize}
\end{thm}

Rather than simply integrate to get $n^i_{\cC/B} = \int_B n^i d \chi$, we can extract
more refined information by taking Chern-Schwarz-Macpherson classes $c^{SM}_*$.  
We recall that $c^{SM}_*$ is the unique map from constructible functions 
to homology which commutes with pushforward and satisfies the normalization
$c^{SM}_*(1_X) = c(TX) \cap [X]$ when $X$ is smooth projective; its existence was
conjectured by Deligne and Grothendieck and established by Macpherson \cite{Mph}. 
Now taking $\mathbf{n}^i := c^{SM}_*(n^i) \in H_*(B)$, we see by Macpherson's theorem 
that 

\[
\sum_{n=0}^\infty \pi^{[n]}_* (c(T \cC^{[n]}) \cap [\cC^{[n]}])  q^{n+1-g} = 
\sum_{i=0}^\infty  \mathbf{n}^i_{\cC/B} \left(\frac{q}{(1-q)^2}\right)^{i+1-g}.
\]

Of course we are also free to take this as the definition of the $\mathbf{n}^i_{\cC/B}$
and conclude from Macpherson's theorem that $\mathbf{n}^i_{\cC/B}$ vanishes if
$i$ is greater than the maximum cogenus of any curve in the family. 
In good cases, the constructible function $n^i$ and the class $\mathbf{n}^i$ carry singularity-theoretic
meaning: 

\begin{thm} \label{thm:euob} 
	Let $\cC/B$ be a family of reduced plane curves with all $\cC_B^{[n]}$ nonsingular.  
	Let $B^i$ be the locus of curves of cogenus $i$, and let $B^i_+$ be the sublocus of curves smooth 
	away from $i$ nodes.  Assume $B^j \subset \overline{B^i_+}$ for all $j \ge i$.  Then
	\[\mathrm{mult}_b (\overline{B^i_+}) = n^i(b) = \mathrm{Eu}_b (\overline{B^i_+}).\] 
	Here
	$\mathrm{Eu}$ is the local Euler obstruction.  Moreover,
	$\mathbf{n}^i_{\cC/B}$ is the Chern-Mather class of $\overline{B^i_+}$. 
\end{thm}
\begin{proof}
The first equality was asserted in \cite{Sh} in the case of a 
locally versal family, but in fact the same argument applies in the above generality.  

For the second, it is shown in \cite{MiSh2} that for {\em any} proper map $f: X \to Y$ of algebraic varieties,
the varieties $V^\alpha$ which appear in the expansion
$f_* 1 = \sum c_\alpha \cdot \mathrm{Eu}(V^\alpha)$
are all components of the higher discriminants (introduced in \cite{MiSh2}) of the morphism $f$. 
It was shown in \cite{Sh} that the higher discriminants of the map $\cC_B^{[n]} \to B$ are precisely 
the loci $\overline{B^i}$, which by assumption agree with $\overline{B^i_+}$. 
It remains to determine the coefficients $c_\alpha$; this may be done at the general point of each
$\overline{B^i_+}$, i.e., in $B^i_+$.  Now we can either make a calculation for nodal curves, 
or alternatively observe that since the loci $B^i_+$ are all immersed in $B$ their Euler obstructions
agree with their multiplicities. 

 The identification of $\mathbf{n}^i$ with the Chern-Mather class of $\overline{B^i_+}$ now follows from
 Macpherson's construction of the functorial Chern class. 
\end{proof}

\begin{rem}
	The assumption on genericity of nodal curves in 
	Theorem \ref{thm:euob} holds for a locally versal family by \cite{DH,T}, 
	for the general $\proj^\delta \subset |L|$ when $L$ is $\delta$-very-ample by \cite{KST},
	for the general $\proj^\delta \subset |L|$ when $L$ is irreducible on a general K3 surface 
	(the genericity of nodal curves in maximal cogenus is \cite{Che}), 
	and it was known classically 
	for the general $\proj \subset |\mathcal{O}_{\proj^2}(d)|$ containing no nonreduced curves. 
\end{rem}

\begin{rem}
    Aluffi \cite{Alu} has shown that the multiplicity and and Euler obstruction 
    of the discriminant of cubic curves on $\proj^2$ differ 
    at a triple line.  The argument above fails here because the total space of the
    restriction of the universal family to a one dimensional disc passing through
    this point is necessarily singular.  Aluffi has also managed to extract enumerative
    information about curves with singularities more complicated than nodes from
    the Chern-Mather and Chern-Schwarz-Macpherson classes of discriminants; perhaps
    the same can be done with the higher Severi strata. 
\end{rem}

The map $c_*^{SM}$ admits a refinement due to Brasselet, Sch\"urmann, and
Yokura.  We denote it $X_{-y}^{BSY}:\mhm(\cdot) \to \mathrm{H}_*(\cdot)[y]$. It 
commutes with pushforward and obeys the normalization 
$X_{-y}^{BSY}(\Q_M) = X_{-y}(TM) \cap [M]$ for $M$ proper smooth.  Thus 
we may apply their functor to the $\mathcal{N}^i_{\cC/B}$ and conclude 
that there are $\mathbf{N}^i_{\cC/B}(y) = X_{-y}^{BSY}(\mathcal{N}^i)
\in \mathrm{H}_*(B)[y]$
such that $\mathbf{N}^i_{\cC/B}(y)=0$ for $i$ greater than the arithmetic
genus of the curves (in a family of integral curves with nonsingular relative Hilbert schemes), and:

\[
\sum_{n=0}^\infty \pi^{[n]}_* (X_{-y}(T \cC^{[n]}) \cap [\cC^{[n]}])  q^{n+1-g} = 
\sum_{i=0}^\infty  \mathbf{N}^i_{\cC/B}(y) \left(\frac{q}{(1-q)(1-qy)} \right)^{i+1-g}.
\]

In the case of interest when $B = \P^\delta$, we denote by abuse 
$\mathbf{N}^i_{\cC/B}(y,H) \in \mathrm{H}^*(\P^\delta)$ the Poincar\'e dual
class of $\mathbf{N}^i_{\cC/B}(y)$.  Since  $\mathbf{N}^i_{\cC/B}(1)$
is, in good cases, the Chern-Mather class of the (codimension i) Severi variety of
cogenus i curves, we might expect:

\begin{conj}\label{cochyfam:stronger}  Let $L$ be a line bundle on a surface $S$, and let
$\CP^\delta \subset |L|$ be a linear subsystem of reduced curves over which the relative
Hilbert schemes $\cC^{[n]} \to \CP^\delta$ are nonsingular for all $n \ge 0$. 
Then $\mathbf{N}^i_{\cC/\CP^\delta}(y,H)$ is a polynomial of minimal degree $i$ in $H$,
and in particular vanishes for $i > \delta$. 
\end{conj}

We will see shortly that this conjecture is in fact equivalent to 
\conjref{cochyfam}.  First we 
recall how the Hirzebruch genera of the relative Hilbert schemes
may be computed.

\subsection{Genera of relative Hilbert schemes} \label{sec:genera}

Following Hirzebruch \cite{Hi}, we take a (normalized) genus to 
mean a
natural transformation of contravariant
functors $\Phi:K^0(\cdot) \to \mathrm{H}^*(\cdot,\Lambda)$ (where $\Lambda$ is a commutative ring) 
such that
\begin{itemize}
  \item For the trivial bundle $\C$, we have $\Phi(\C) = 1$.
  \item Sums go to products: $\Phi(E \oplus F) = \Phi(E) \Phi(F)$.
  \item There is a power series $f_\Phi \in 1 + z\Lambda\lbr z\rbr $ 
  such that, for a line bundle $L$, we have $\Phi(L) = f_\Phi(c_1(L))$.
\end{itemize}
In the remainder of the paper we will be concerned only with 
the Hirzebruch genus $\Phi = X_{-y}$, for which 
$\Lambda = \Q\lbr y\rbr $ and $f(z) = \frac{z(1-ye^{-z(1-y)})}{(1-e^{-z(1-y)})}$.
In any case fix some $\Lambda, \Phi$.  We write $\phi(X):=\Phi(TX)$. 

Let $S$ be a surface, $L$ a line bundle on it, $\CP^\delta \subset |L|$ 
some linear system, $H = \oo_{\CP^\delta}(1)$. 
Let $S^{[n]}$ be the Hilbert scheme of $n$ points on $S$, 
and let $Z_n(S)\subset S\times S^{[n]}$ be the universal family, 
with the projections $q:Z_n(S)\to S$, $p:Z_n(S)\to S^{[n]}$. 
Let $L^{[n]}:=p_*q^*L$.  This is a vector bundle of rank 
$n$ on $S^{[n]}$ with fibre $\mathrm{H}^0(Z,L|_Z)$ over $Z\in S^{[n]}$. 
Let $\pi:\cC_{\CP^\delta}\to \CP^\delta$ be the universal curve over $\CP^\delta$ 
and denote by $\pi^{[n]}:\cC_{\CP^\delta}^{[n]} \to \CP^\delta$ 
the relative Hilbert scheme of points.  
The relative Hilbert scheme always has the expected dimension $\delta + n$
\cite{AIK},
and is the scheme theoretic zero locus of a tautological section 
of $L^{[n]}\boxtimes H$; when $\cC_{\CP^\delta}^{[n]}$ 
is nonsingular this section is transverse.

As a book-keeping device, let $e^x$ denote a trivial 
line bundle with nontrivial $\C^*$ action giving
equivariant first Chern class $x$, i.e., 
$\Phi(e^x)=f_\Phi(x)$.  



\begin{defn}
  Let \[D_n^{S,L,\Phi}(x):=  \int_{S^{[n]}} 
    \Phi(TS^{[n]})  \frac{c_n(L^{[n]} \otimes e^x)}
  {\Phi(L^{[n]} \otimes e^x)}  \in \Lambda\lbr x\rbr. \]
\end{defn}

\begin{prop} \label{prop:hilbgenus}
  Assume $\cC^{[n]}_{\CP^\delta}$ is nonsingular.  Then 
  \[\pi^{[n]}_* (\Phi(T\cC^{[n]}_{\CP^\delta}) \cap [\cC^{[n]}]) = 
    f_\Phi(H)^{\delta+1}
  D_n^{S,L,\Phi}(H) \cap [\P^\delta] .\]
  In particular,
  \[\phi(\cC^{[n]}_{\CP^\delta}) = \Res_{x=0} \left( \frac{f_\Phi(x)}{x}\right)^{\delta+1}
  D_n^{S,L,\Phi}(x) .\]
\end{prop}

\begin{proof}
  Denote $q:S^{[n]}\to pt$.  Then 
\begin{eqnarray*}
  & & \pi^{[n]}_* (\Phi(T\cC^{[n]}_{\CP^\delta}) \cap [\cC^{[n]}]) \\
& = 
    & (q \times 1_{\P^\delta})_* \left( c_n(L^{[n]} \boxtimes H) 
    \frac{\Phi(TS^{[n]}) \Phi(T\CP^{[\delta]})}
    {\Phi(L^{[n]} \boxtimes H)} \cap [S^{[n]}] \right) \cap [\P^\delta] \\
  & = &  \Phi(H)^{\delta+1} q_*\left[  \left.
    \Phi(TS^{[n]})  \frac{c_n(L^{[n]} \otimes e^x)}
  {\Phi(L^{[n]} \otimes e^x)} \right]\right|_{x=H} \cap [\P^\delta] .
\end{eqnarray*}
\end{proof}

 \begin{rem}
   If $h^0(L)>\delta$,  the formula makes sense without requiring smoothness,
  if we view it as a virtual contribution. The description of 
  $\cC_\delta^{[n]}$ as zero locus of a section of 
  $L^{[n]}\boxtimes H$ gives it a virtual fundamental class 
  and a virtual tangent bundle (see e.g. \cite{FG}). Thus 
  independent of the singularities of $\cC_\delta^{[n]}$, what 
  is computed here is $\phi(\cC^{[n]})$ 
  with this virtual structure. 
  Without any assumption on $L$, we can view  the second equality in 
 Proposition 24 as a definition of $\phi(\cC^{[n]})$.
\end{rem}

The $D^{S,L,\Phi}_n$ enjoy a certain multiplicativity.  Introduce the series
\begin{equation}
  D^{S,L,\Phi} = \sum D^{S,L,\Phi}_n q^n.
\end{equation}

We  denote by $[S,L]$ the algebraic cobordism class of the 
pair of the surface $S$ and the line bundle $L$.
For us\footnote{
For a geometric account of algebraic cobordism of varieties with bundles see \cite{LP}, we however
do not require any results of this theory and only use ``cobordism class'' as a convenient shorthand.}
this is the equivalence class of pairs $(S,L)$, 
where two such $(S_1,L_1)$, $(S_2,L_2)$ are equivalent 
if the numbers $L_i^2$, $L_i K_{S_i}$, $K_{S_i}^2$, $c_2(S_i)$ 
coincide for $i=1,2$. 
These form a group where in particular $[S_1,L_1]+[S_2,L_2]=[S_1\sqcup S_2,L]$
 with $L$ the line bundle which is $L_1$ on $S_1$ and $L_2$ on $S_2$.

\begin{prop} \label{DLprod}
  The map $[S,L] \mapsto D^{S,L,\Phi}$ is a homomorphism from
  the cobordism group of surfaces with bundles to 
  the multiplicative group of invertible power series in $q$ with coefficients in $\Lambda$.  
  In particular there exist
  universal power series $D_1,D_2,D_3,D_4\in \Lambda\lbr q\rbr $ such that 
  $$D^{S,L,\Phi}=D_1^{L^2} D_2^{L K_S} D_3^{K_S^2} D_4^{c_2(S)}.$$
\end{prop}
\begin{proof}
  The $D^{S,L,\Phi}_n(x)$ are defined by a genus applied to 
  $L^{[n]}$ and $T_{S^{[n]}}$, so the first
  statement follows from the arguments in \cite{EGL}.  

 For the second, let $(S_i,L_i)$, $i=1,\ldots,4$ be chosen such that the 
$4$ vectors $(L_i^2, L_i K_{S_i}, K_{S_i}^2, c_2(S_i))$ 
are linearly independent.  Then 
for any $(S,L)$ we can write $[S,L]=a_1[S_1,L_1]+a_2[S_2,L_2]+
a_3[S_3,L_3]+a_4[S_4,L_4]$, with $a_i\in \Q$ and thus 
$$D^{S,L,\Phi}=\prod_{i=1}^4(D^{S,L_i,\Phi})^{a_i}.$$ 
For simplicity we write $[S_i,L_i]$ also for the corresponding vector.
Let $e_1=(1,0,0,0),\ e_2=(0,1,0,0),\ e_3=(0,0,1,0),\ e_4=(0,0,0,1)$. Choosing 
$(a_{i,j})_{i,j=1}^4$ in $\Q$ with $\sum_{j}a_{i,j}[S_j,L_j]=e_i$ for all $i$, we put 
$D_i=\prod_{j=1}^4(D^{S,L_i,\Phi})^{a_{i,j}}$. As $[S,L]=L^2e_1+LK_S e_2+K_S^2 e_3+e(S)e_4$, this gives  
$D^{S,L,\Phi}=D_1^{L^2} D_2^{L K_S} D_3^{K_S^2} D_4^{c_2(S)}.$
\end{proof}

The multiplicativity of Proposition \ref{DLprod} allows  
$D^{S,L,\Phi}_n(x)$ (and therefore also
 $\phi(\cC^{[n]})$ and the BPS invariants) 
to be computed by localization in the following standard way.

By the above argument to compute $D^{S,L,\Phi}$ for any $(S,L)$ it is enough to 
compute them for pairs $(L_i,S_i)$ whose corresponding vectors are linearly independent .
We can choose the $(S_i,L_i)$ as $(\P^2,\oo)$, $(\P^1\times \P^1,\oo)$,  
$(\P^1\times \P^1,\oo(1,0))$, $(\P^1\times \P^1,\oo(1,-1))$.
In this case $S=S_i$ is a toric surface, i.e. it has an action by 
a torus $(\C^*)^2$ with finitely many fixed points, and $L=L_i$ has 
an natural equivariant lifting. The action of $(\C^*)^2$ on $S$ 
induces in a natural way an action 
on $S^{[n]}$, and the equivariant lifting of $L$ induces an 
equivariant lifting of $L^{[n]}$.  Thus we can apply equivariant
 localization to compute $D^{S,L,\Phi}_n(x)$, in terms of the weights 
of the action on the fibres of $T_{S^{[n]}}$ and $L^{[n]}$ at the 
fixed points. The fixed points are parametrized by tuples 
of Young diagrams and the weights of the action can be expressed explicitly in 
terms of this data.
For more details in a slightly different situation see e.g. \cite{EG,NY, CO}.

From now on we specialize to $\Phi = X_{-y}$ and abbreviate
$D^{S,L}_n(y,x):=D^{S,L,X_{-y}}_n(x)$, $D^{S,L}(y,x,q):=\sum_{n\ge 0} D^{S,L}_n(y,x)q^n$.
In this case  a computer calculation yields the 
$D^{S,L}_n(y,x)$ for $n\le 10$ and modulo $x^{14}$. 
The $\chi_{-y}(\cC_{\CP^\delta}^{[n]})$ are computed from this by \propref{prop:hilbgenus}.

\subsection{A reformulation of the conjectures, and evidence.}

Let $t_i$ be the Chern
roots of the tangent bundle of the Hilbert scheme, and let $l_i$ be
the Chern roots of the bundle $L^{[i]}$.  In the previous subsection we introduced the series
\begin{equation}\label{DSL}
 D^{S,L}(y,x,q):=
 \sum q^n \int_{S^{[n]}}
  \prod_{i=1}^{2n} \frac{t_i(1-ye^{-t_i(1-y)})}{(1-e^{-t_i(1-y)})}
  \prod_{j=1}^n \frac{(1-e^{-(l_j+x)(1-y)})}{(1-ye^{-(l_j+x)(1-y)})} \in \Q[y]\lbr x\rbr \lbr q\rbr .
\end{equation}
For convenience we write $\que  = q/((1-q)(1-qy))$.\footnote{
We record here that the compositional inverse is given by what are 
called the Narayana numbers, 
\[q(Q) = \sum_{n=1}^\infty \sum_{k=1}^n Q^n  y^{k-1} \frac{(-1)^{n-1}}{n} {n \choose k} {n \choose k-1}\]
which specializes to the following formulas involving Catalan numbers,
\begin{eqnarray*}
q(Q)|_{y=1} & = & \sum_{n=1}^\infty Q^n \frac{(-1)^{n-1}}{n+1} {2n \choose n} \\
q(Q)|_{y=0} & = & \frac{Q}{1+Q}  \\
q(Q)|_{y=-1} & = & \sum_{n=0}^\infty Q^{2n+1} \frac{(-1)^n}{2n+1}{4n \choose 2n} .
\end{eqnarray*}
Note also that $Q(q) \in q\Z\lbr q,qy\rbr $ and $q(Q) \in Q\Z\lbr Q,Qy\rbr .$
}

By Proposition \ref{prop:hilbgenus} we have

\begin{equation}\label{Dpow}
  \sum_{i=0}^\infty  \mathbf{N}^i_{\cC/\CP^\delta}(y,H)
  \que^{i} = \left(\frac{q}{Q}\right)^{1-g} 
  \left(\frac{H(1-y e^{-H(1-y)})}{1-e^{-H(1-y)}}\right)^{\delta+1}
  D^{S,L}(y,H,q).
\end{equation}

\begin{conj}\label{DLconj} For any surface $S$ and line bundle $L$, we have
\[ (q/Q)^{1-g(L)} D^{S,L}(y,x,q)  \in \Q[y]\lbr x, xQ\rbr  .\]
\end{conj}

\begin{prop}
  Conjectures \ref{cochyfam}, \ref{cochyfam:stronger}, and \ref{DLconj} are equivalent. 
\end{prop}
\begin{proof}
  Note $f(x):=X_{-y}(e^x) \in 1 + x\Q[y]\lbr x\rbr $ is invertible.  From Equation \ref{Dpow},
  we see that
  \conjref{DLconj} implies 
  \conjref{cochyfam} and \conjref{cochyfam:stronger}.

  Assume \conjref{cochyfam}.  Now 
  consider some fixed linear system $\P^\delta \subset |L|$ on some surface $S$ 
  such that all the relative Hilbert schemes $\cC^{[n]} \to \P^\delta$ are nonsingular.
  \conjref{cochyfam:stronger} amounts to the statement 
  \[(Q/q)^{g-1} f(x)^{\delta + 1} D^{S,L}(y,x,q) \in \Q\lbr x,xQ\rbr  +O(x^{\delta+1}).\]
  This obviously holds at $\delta = 0$; let us prove it holds at $\delta$ by induction.
  If we know this statement holds for some
  $\delta = r-1$ and wish to check it for $\delta = r$, since 
  $f(x) \in \Q\lbr x\rbr $ we already know the statement modulo $x^r$.
  So we need only check 
  \[\deg_Q \Coeff_{x^r} 
  (Q/q)^{g-1} f_\Phi(x)^{r + 1} D^{S,L}(y,x,q) \le r .\]
  But this is precisely the assertion of \conjref{cochyfam} for $\P^r$.  
  And since all the relative Hilbert schemes will also be smooth over
  a general $\P^r \subset \P^\delta$ for any $r \le \delta$, the hypothesis
  of \conjref{cochyfam} is satisfied and we may deduce \conjref{cochyfam:stronger}.

  Finally, Conjecture \ref{cochyfam:stronger} asserts that Conjecture \ref{DLconj} 
  holds modulo $x^{\delta+1}$.  We 
  have an expression $D^{S,L} = D_1^{L^2} D_2^{L.K_S} D_3^{K_S^2} D_4^{c_2(S)}$
  where the $D_i$ are power series starting with $1$. Thus we compare $D^{S,L}$ for various surfaces
  and line bundles 
  to conclude the statement of Conjecture \ref{DLconj}
  for the series $D_i$, modulo some $x^k$.  Taking $k \to \infty$ by choosing 
  increasingly ample line bundles recovers the statement for $D^{S,L}$. 
\end{proof}

\begin{rem}
  The above argument implies in particular that 
  Conjecture \ref{DLconj} holds in the Euler charateristic
  limit $y=1$.  From this it follows formally that for {\em any}
  $\P^\delta \subset |L|$, with no assumptions on the reducedness
  or irreducibility of the curves that appear or on the smoothness
  of the relative Hilbert schemes, there are integers $n^i_{\cC/\P^\delta}$ 
  such that \[\sum_{n=0}^\infty q^{n+1-g} \int_{[\cC^{[n]}]^{vir}} 
  c_{top}(T^{vir} \cC^{[n]}) = 
  \sum_{i=0}^{\delta} n^i_{\cC/\P^\delta} \left(\frac{q}{(1-q)^2}\right)^{i+1-g}. \]
\end{rem}

Recall that $D^{S,L}$ can be expressed in four universal power series,
\[D^{S,L} = D_1^{L^2} D_2^{K_S.L} D_3^{K_S^2} D_4^{c_2(S)}.\]  
To avoid writing $(Q/q)^{g-1}$ we adjust these series slightly.

\begin{defn}
We write $\widetilde{D}^{S,L} := (Q/q)^{g-1} D^{S,L}$.  We also take 
$\widetilde{D}_1 := (Q/q)^{1/2} D_1$,
$\widetilde{D}_2:= (Q/q)^{1/2} D_2$, and $\widetilde{D}_3 = D_3$,
$\widetilde{D}_4 = D_4$ in order that 
$\widetilde{D}^{S,L} = \widetilde{D}_1^{L^2} \widetilde{D}_2^{K_S.L} \widetilde{D}_3^{K_S^2} 
\widetilde{D}_4^{c_2(S)}$.
\end{defn}

We have $\widetilde{D}^{S,L} \in 1 + (y,x,Q)\Q[y]\lbr x,Q\rbr $ for all 
$S,L$, hence the same is true for the $\widetilde{D}_i$. 
Similarly, \conjref{DLconj} is equivalent to the assertion that
$\widetilde{D}_i \in \Q[y]\lbr x, xQ\rbr $ for all $i$.

\begin{thm} \label{thm:vanishing}
  $\widetilde{D}_1, \widetilde{D}_4 \in \Q[y]\lbr x,xQ\rbr $. 
\end{thm} 
\begin{proof}
  Let $(A, L)$ be a primitively polarized abelian surface of Picard rank 1.  If $L^2 = 2k+2$ then
  $\dim |L| = k$ and the curves in $|L|$ have arithmetic genus $k + 2$.  Note such $(A, L)$ exist for all $k$. 
  By \cite{Muk} the relative Hilbert schemes are smooth, and so from Proposition \ref{prop:mhmbps} 
  we find that the $N^i$
  vanish beyond the arithmetic genus.  By this vanishing and the formula (\ref{Dpow}) extracting the $N^i$ from 
  $D^{S,L}$: 
  \[\deg_Q \Coeff_{x^k} 
  f_\Phi(x)^{k + 1} \widetilde{D}_1^{2k+2} \le k+2  .\] 
   We write $\xi(x,Q) = f_\Phi(x) \widetilde{D}_1^2 \in 1 + (x,Q)\Q[y]\lbr x,Q\rbr $. We want to show
   $\widetilde{D}_1 \in \Q[y]\lbr x,xQ\rbr $; since this evidently holds for $f_\Phi(x)$ and we may 
   take roots of power series starting with $1$, 
 it suffices to show this for $\xi$.  So we have 
   $\deg_Q \Coeff_{x^k} \xi(x,Q)^{k+1} \le k+2$.

   The following argument is completely formal and does not involve the geometric meaning of $\xi$.    
   We write $d_Q(k):= \deg_Q \Coeff_{x^k} \xi(x,Q)^{k+1}$. 
   Let $k_1 = \min \{k\,|\,d_Q(k) > k \}$, assuming this set is nonempty. 
   Then $d_Q(k_1) = \deg_Q \Coeff_{x^{k_1}} \xi(x,Q)^{{k_1}+1}$,
   since
   no lower (in $x$) degree term can contribute such a high power of $Q$.  There are two cases, 
   $d_Q(k_1) = k_1 + 2$ or $d_Q(k_1) = k_1 + 1$.  In the first case, consider 
   $\Coeff_{x^{2k_1}} \xi(x,Q)^{2k_1 +1}$.  There will be a contribution from products
   of two terms of the form $Q^{k_1+2} x^{k_1}$, which gives the highest possible power of $Q$ and
   thus cannot be canceled.  But then $\deg_Q \Coeff_{x^{2k_1}} \xi(x,Q)^{2k_1 +1} = 2k_1+4$, which
   is a contradiction.  In the second case, consider $\Coeff_{x^{3k_1}} \xi(x,Q)^{3k_1 +1}$.  In order
   that the degree $3k_1 + 3$ contribution from products of three terms $Q^{k_1+1} x^{k_1}$ be cancelled,
   there must be some $h + h' = 3k_1$ with $d_Q(h) = h +2$ and $d_Q(h') > h'$. By minimality
   of $k_1$, we have $h' > k_1$ hence $h < 2k_1$.   Let $k_2 = \min \{k\,|\,d_Q(k) > k+1 \} \le h < 2k_1$.
   Finally consider $\Coeff_{x^{k_1 + k_2}} \xi(x,Q)^{k_1 + k_2 +1}$.  There is a contribution from
   products of terms the form $x^{k_1} Q^{k_1+1}$ and $x^{k_2} Q^{k_2+2}$; since $k_2 < 2k_1$ this
   contribution cannot be cancelled.  This is a contradiction.  So finally we must have 
   $d_Q(k) \le k$ for all $k$, hence $\xi(x,Q) \in \C[y]\lbr x,xQ\rbr $, hence the same holds for $\widetilde{D}_1$. 
   
   Now let $(K, L)$ be a primitively polarized K3 surface of Picard rank 1.  If $L^2 = 2g-2$ then
   $\dim |L| = g$ and the curves in $|L|$ have genus $g$.  Such $(K, L)$ exist for all $g$.  By vanishing
   of the $N^i$ beyond the arithmetic genus we have 
     \[\deg_Q \Coeff_{x^g} 
  f_\Phi(x)^{g + 1} \widetilde{D}_1^{2g-2} D_4^{24} \le g . \]
    Since we know $f_\Phi(x), \widetilde{D}_1 \in \Q[y]\lbr x,xQ\rbr $, we may conclude the same for $D_4$. 
\end{proof}

\begin{cor}
   Conjectures \ref{cochyfam}, \ref{cochyfam:stronger}, and \ref{DLconj} hold for surfaces with numerically trivial canonical class. 
\end{cor}

\begin{rem}
   Note the slightly curious nature of the proof of the theorem and 
   corollary: for geometric reasons, namely smoothness of the relative
   Hilbert schemes and the near equality of the genus and dimension of the linear system for certain line bundles on 
   K3 and abelian surfaces, we know the conjecture for {\em complete} linear systems on K3 surfaces and something
   close for abelian surfaces.  Then by leveraging the universality of the expressions, and the existence of K3 and abelian
   surfaces of all genera, we can conclude the result also for {\em not-necessarily-complete} linear systems. 
   
   This sort of approach was suggested to the authors by Pandharipande \cite{Pan}, who further suggested that the 
   other power series may be similarly constrained by finding enough other surfaces with nontrivial canonical class
   but for which nonetheless the genus and dimension of linear systems are close.  However to use these (or any)
   surfaces for the present purposes, one must establish smoothness of the relative Hilbert schemes, which we 
   do not know how to do.
\end{rem}

Using the localization calculation described in Section \ref{sec:genera},
we can give evidence for \conjref{DLconj} for arbitrary surfaces. 
 
\begin{prop} \conjref{DLconj} holds modulo $q^{11}$ and $x^{14}$.
Therefore, 
if  $0\le \delta\le 13$ and $\P^\delta \subset |L|$ is a linear system
over which the relative Hilbert scheme is smooth, 
there exist polynomials $N^i_{\cC/\P^\delta} (y)
\in \Z[y]$,
where $i = 0,\ldots, \delta$, 
such that 
$$ \sum_{n\ge 0} \chi_{-y}(\cC^{[n]})q^{n+1-g}\equiv
\sum_{l=0}^\delta N^i_{\cC/\P^\delta} (y)
Q^{l+1-g}\mod O(q^{11 + 1 - g}), $$
and furthermore these polynomials are explictly computed.
If moreover $g < 11$ and all curves are irreducible, 
then the equality is established to all orders. 
\end{prop}

For example, Conjecture \ref{cochyfam} holds for a general 
$\P^4$ in $|\mathcal{O}_{\P^2}(6)|$. 

\vspace{2mm} 

The relation between the various $n, \mathbf{n}, N, \mathbf{N}$ and the series $D^{S,L}$ has
in this section always been contingent on the smoothness of the relative Hilbert schemes over 
the appropriate linear subsystem $\proj^\delta \subset |L|$.  To avoid continually making this
hypothesis, we introduce the following:

\begin{defn}\label{Ndef}
For a surface $S$ and a line bundle $L$, we define $\mathbf{N}^i_{\delta,[S,L]}$, $N^i_{\delta,[S,L]}$,  by the formulas
\[
  \sum_{i=0}^\infty  \mathbf{N}^i_{\delta,[S,L]}
  \que^{i} =   \left(\frac{H(1-y e^{-H(1-y)})}{1-e^{-H(1-y)}}\right)^{\delta+1}
  \widetilde{D}^{S,L}(y,H,q),
\]
\[
  \sum_{i=0}^\infty N^i_{\delta,[S,L]}
  \que^{i} =   \Res_{x=0}\left[\left(\frac{1-y e^{-x(1-y)}}{1-e^{-x(1-y)}}\right)^{\delta+1}
  \widetilde{D}^{S,L}(y,x,q)\right],
\]
and similarly for the specializations $n, \mathbf{n}$. 
\end{defn} 

By comparison with Equation \ref{Dpow}, we see that for a linear system
$\CP^\delta \subset |L|$ containing only reduced curves, and whose
relative Hilbert schemes are smooth,  $\mathbf{N}^i_{\CP^\delta \subset |L|} = 
\mathbf{N}^i_{\delta,[S,L]}$.

\section{The term of the deepest stratum} \label{sec:deepest}


For a linear system $\P^\delta \subset |L|$, the 
numbers $n^\delta_{\cC/\CP^\delta}$ have the clearest 
enumerative significance, counting the number of $\delta$-nodal curves 
in the linear system.  
Thus we might also hope that the $N^\delta_{\cC/\CP^\delta}(y)$ 
have an enumerative meaning refining this.  In any case, 
assuming \conjref{DLconj}, the  $N^\delta_{\delta,[S,L]}(y)$ 
are the easiest to compute and their generating function is multiplicative.

We have been in the meantime in able to prove the main conjecture of this section (\conjref{Gconj}) in case $K_S$ is numerically trivial. 
The proof appears in \cite{GSh}, and depends on several results from this paper, Theorem \ref{thm:vanishing} in particular.

\begin{prop}
Assume \conjref{DLconj}, or $K_S = 0$.  
View $\widetilde{D}^{S,L}$ 
as an element of $\Q[y]\lbr x,s\rbr $ with $s=xQ$.
Then $$\sum_{\delta\ge 0} N^\delta_{\delta, [S,L]} s^\delta = \widetilde{D}^{S,L}(y,x=0,s).$$
\end{prop}

\begin{cor} \label{produc} Assume \conjref{DLconj}, or $K_S = 0$.
Then there exist series $A_i \in \Q[y]\lbr s\rbr $ such that
$$\sum_{\delta\ge 0} N^\delta_{\delta,[S,L]} (y) s^\delta=A_1^{L^2}A_2^{LK_S}A_3^{K_S^2}A_4^{c_2(S)}.$$
\end{cor}
\begin{proof}
Viewing $\widetilde{D}_i \in \Q[y]\lbr x,s\rbr $ where $s=xQ$, take $A_i := \widetilde{D}_i|_{x=0}$. 
\end{proof}

In the unrefined ($y = 1$) setting, more explicit formulas were expressed after 
substituting for $s$ a certain quasimodular form.  Specifically, in \cite[Conj. 2.4]{G}, the
following expansion was proposed: 

\begin{equation} \label{Gform-unref}
\sum_{\delta\ge 0}
n^\delta_{\delta,[S,L]} \cdot (DG_2)^\delta=\frac{(DG_2/q)^{\chi(L)}B_1^{K_S^2}B_2^{LK_S}}
{(\Delta \cdot D DG_2/q^2)^{\chi(\oo_S)/2}}.
\end{equation}
Here $G_2$ is the Eisenstein series, $\Delta$ is the discriminant: 
\begin{eqnarray*}
\Delta(q) & = & q \prod_{n = 1}^{\infty} (1-q^k)^{24}, \\
G_2(q) & = & \frac{1}{24} + \sum_{m =1}^\infty q^m \sum_{d|m} d.
\end{eqnarray*}
The series $B_1, B_2 \in 1 + q\Q\lbr q\rbr $ are not known explicitly, 
although their first several coefficients may be computed and are given in \cite{G}. 
We also write $D = q \frac{d}{dq}$. 
The above formula is by now a theorem, since the existence of universal formulas has been established
\cite{Tz, KST}, and the case of the K3 surface (where $K_S = 0$ hence the $B_i$ do not appear)
was solved explicitly \cite{BL}.

\begin{nota} We write
$\overline N^\delta_{\delta,[S,L]}:=N^\delta_{\delta,[S,L]}/y^{\delta}$.
\end{nota}

Now we give a conjectural refinement of Equation \ref{Gform-unref}.  The series
$\Delta, DG_2$ are refined as follows:\footnote{
These functions are related to certain Jacobi forms.  Let $q=e^{2\pi i \tau}$, $y=e^{2\pi i z}$ 
\begin{enumerate}
\item $\widetilde\Delta(y,q)=\phi_{10,1}(\tau,z)/(y^{1/2}-y^{-1/2})^2$.
Here $\phi_{10,1}(\tau,z)=\eta(\tau)^{18}\theta(\tau,z)^2$ is up to normalization the unique Jacobi cusp form on $Sl_2(\Z)$ of weight 10 and index $1$.
\item We can write 
\begin{align*}
(y-2+y^{-1})\widetilde{DG}_2&=\sum_{m= 1}^\infty q^{m} \sum_{d | m} \frac{m}{d} 
(y^d -2+ y^{-d})\\
&=-2(G_2(\tau)+1/24)+\sum_{d,e>0} e(y^d-y^{-d})q^{de}\\
&=-\frac{1}{2}D \log\left(\frac{\phi_{10,1}(\tau,z)}{\Delta(\tau)}\right)=-\frac{1}{2}D \log\left(\phi_{-2,1}(\tau,z)\right)\!.
\end{align*}
Here $\Delta(\tau)$ is the discriminant function and $\phi_{-2,1}=\phi_{10,1}/\Delta$ is the up to normalization unique weak Jacobi cusp form  of weight $-2$ and index $1$ on $Sl_2(\Z)$.
\end{enumerate}
}
\begin{align*}
\widetilde\Delta(y,q)&:= q \prod_{n=1}^{\infty}(1-q^n)^{20}(1-yq^n)^{2}
(1-y^{-1}q^n)^2,\\
\widetilde{DG}_2&:=\sum_{m= 1}^\infty m q^{m} \sum_{d | m} \frac{[d]_y^2}{d} .
\end{align*}

\begin{conj}
\label{Gconj}
There exist universal power series 
$B_1(y,q)$, $B_2(y,q)$ in $\Q[y,y^{-1}]\lbr q\rbr $, such that 
\begin{equation}
\label{Gform}
\sum_{\delta\ge 0}
\overline{N}^\delta_{\delta, [S,L]}(y) (\widetilde{DG}_2)^\delta=\frac{(\widetilde{DG}_2/q)^{\chi(L)}B_1(y,q)^{K_S^2}B_2(y,q)^{LK_S}}
{\big(\widetilde\Delta(y,q)\, D \widetilde{DG}_2/q^2)^{\chi(\oo_S)/2}}.
\end{equation}
Here, to make the change of variables, all functions are viewed as elements of 
$\Q[y,y^{-1}]\lbr q\rbr $.
\end{conj}
This conjecture is again checked modulo  $q^{11}$; and  we get that
{\small\begin{align*}
&B_1(y,q)=1 - q - ((y^2 +3y + 1)/y)q^2 
+ ((y^4 + 10y^3 + 17y^2 + 
10y+ 1)/y^2)q^3  
 - ((18y^4 + 87y^3 \\ &+ 135y^2 + 87y + 18)/y^2)q^4
+ ((12y^6 + 210y^5 + 728y^4 + 1061y^3 + 728y^2 + 210y + 12)/y^3)q^5 \\ &- ((2y^8 + 259y^7 + 2102y^6 + 5952y^5 + 8236y^4 + 5952y^3 + 2102y^2 + 259y + 2)/y^4)q^6 + ((162y^8 \\ &+ 3606y^7  + 19668y^6 + 48317y^5 + 64253y^4 + 48317y^3 + 19668y^2 + 3606y + 162)/y^4)q^7  \\ &
- ((47y^{10} + 3789y^9 + 41999y^8 + 177800y^7 + 392361y^6 + 505678y^5 
 + 392361y^4 + 177800y^3 \\ &+ 41999y^2 + 3789y + 47)/y^5)q^8 
+ ((5y^{12} + 2416y^{11}  + 60202y^{10} + 445989y^9 + 1576410y^8 \\ &+ 3197831y^7 
+ 4018919y^6 + 3197831y^5  + 1576410y^4 + 445989y^3 + 60202y^2 + 2416y + 5)/y^6)q^9\\ & -
((896y^{12} + 58504y^{11} + 793194y^{10} + 4483755y^9 + 13818256y^8 +26192369y^7 + 32243357y^6 \\ &+ 26192369y^5 + 13818256y^4 +4483755y^3 + 793194y^2 + 58504y + 896)/y^6)q^{10}+O(q^{11}),
\\
 &B_2(y,q)=\frac{1}{(1-yq)(1-q/y)}\big(1 + 3q - ((3y^2 + y + 3)/y)q^2 + ((y^4 + 8y^3 + 18y^2 + 8y + 1)/y^2)q^3\\ &- ((13y^4 + 53y^3 + 76y^2 + 53y + 13)/y^2)q^4 + ((7y^6 + 100y^5 + 316y^4 + 455y^3 + 316y^2+ 100y  \\&+ 7)/y^3)q^5- ((y^8 + 112y^7 + 779y^6 +2076y^5 + 2819y^4 + 2076y^3 + 779y^2 + 112y + 1)/y^4)q^6\\& + ((67y^8 + 1243y^7 + 6129y^6 + 14386y^5 + 18870y^4 + 14386y^3 + 6129y^2 + 1243y + 67)/y^4)q^7\\&- ((19y^{10} + 1281y^9 + 12417y^8 + 48879y^7 + 104034y^6 + 132579y^5 +104034y^4 + 48879y^3 \\&+ 12417y^2 + 1281y + 19)/y^5)q^8 + ((2y^{12} + 822y^{11} + 17542y^{10} + 117829y^9 + 393703y^8 \\&+ 775411y^7 + 965540y^6 + 775411y^5 + 393703y^4 + 117829y^3 + 17542y^2 + 822y + 2)/y^6)q^9\\& -((310y^{12} + 17206y^{11} + 207074y^{10} + 1085712y^9 + 3197506y^8 + 5913778y^7 + 7223539y^6 \\&+ 5913778y^5 + 3197506y^4+ 1085712y^3 + 2070742y^2 + 17206y + 310)/y^6)q^{10} +O(q^{11})\big).\\
\end{align*}}

At $y=1$, we recover modulo $q^{11}$ the functions $B_1(q)$, $B_2(q)$ of 
\cite{G}.

\vspace{2mm} As in \cite[Rem. 2.6]{G}, the expansion in $\widetilde{DG}_2$ may be exchanged for
an expansion in $q$ while simultaneously trading a sum over varying numbers of point conditions 
while fixing the line bundle for a sum
over line bundles while fixing the point conditions.   Note the latter form is more natural from the point of
view of the GW/DT/pairs theories, and indeed this is the form in which the K3 is solved in \cite{BL}. 

In detail this procedure is as follows.  For any power series $f \in R\lbr q\rbr $ and $g \in q + q^2 R\lbr q\rbr $, we
may expand $f$ in terms of $g$ by the residue formula:
\[f(q) = \sum_{l = 0}^\infty g(q)^l \Coeff_{q^0}\left[\frac{f(q)Dg(q)}{g(q)^{l+1}}\right].\]
Conjecture \ref{Gconj} asserts that $\overline{N}^\delta_{\delta,[S,L]}$ is the coefficient of 
$\widetilde{DG}_2^\delta$ of a certain expression; taking this coefficient by the above residue formula
gives the equivalent formulation:
\begin{eqnarray*}
\overline{N}^\delta_{\delta, [S,L]}(y) & = &
\Coeff_{q^0}\left[\widetilde{DG}_2^{-\delta-1} D\widetilde{DG}_2
\frac{(\widetilde{DG}_2/q)^{\chi(L)}B_1(y,q)^{K_S^2}B_2(y,q)^{LK_S}}
{\big(\widetilde\Delta(y,q)\, D \widetilde{DG}_2/q^2)^{\chi(\oo_S)/2}}\right] \\
& = & \mathrm{Coeff}_{q^{(L^2 - LK_S)/2}}\left(
\frac{(\widetilde{DG}_2)^{\chi(L)-1-\delta} D \widetilde{DG}_2 B_1(y,q)^{K_S^2}B_2(y,q)^{LK_S}}
{\big(\widetilde\Delta(y,q)\, D \widetilde{DG}_2)^{\chi(\oo_S)/2}} \right)\!.
\end{eqnarray*}

We would now like to collect coefficients of $q$ to write the entire series in the $( \cdot )$ in terms of
the $\overline{N}$.  So we must choose
some values of $\delta, [S, L]$ such that $K^2_S, LK_S, \chi(\oo_S)$ and $k:=\chi(L)-1-\delta$ remain constant, 
but $(L^2-LK_S)/2$ assumes every integer value starting from $k +1 -  \chi(\oo_S)$.   In other words, the cobordism
class of the surface, the number of point conditions $k = \chi(L)-1-\delta$, and 
$LK_S$ are fixed, and $L^2$ varies.  Note that it is not necessarily possible to find a fixed surface $S$
and honest line bundles $L_i$ which realize all the desired values.  This causes no difficulties as the
$\overline{N}^{\delta}_{\delta,[S,L]}$ may be viewed as just functions of the four values $L^2, LK_S, K_S^2, c_2(S)$. 
Making this dependence explicit we write 
$\overline{M}_{k,[S]}((L^2-LK_S)/2, LK_S):=\overline{N}^{\chi(L)-1-k}_{\chi(L)-1-k,[S,L]}$, where the RHS is viewed just
as a function of two integers and is determined by evaluating the LHS on the cobordism class $[S,L]$ 
with the specified invariants.  In terms of the $\overline{M}$, we have: 

\begin{equation} \label{eq:Ms}
\sum_{l = k+1-\chi(\oo_S)}^{\infty}\!\!\!\!\!\!\!\!\! \overline{M}_{k,[S]}(l,LK_S)q^l = \widetilde {DG}_2(y,q)^k 
\frac{ B_1(y,q)^{K_S^2}B_2(y,q)^{LK_S} D \widetilde {DG}_2(y,q)}{(\widetilde \Delta(y,q)D\widetilde {DG}_2(y,q))^{\chi(\oo_S)/2}}.
\end{equation}

We can also express this in a slightly different way, fixing $LK_S$, $k$ and varying $\delta$: write $\overline N_{k,[S]}(\delta,LK_S):=\overline N^{\delta}_{\delta,[S,L]}$ with $k=\chi(L)-1-\delta$. Then
\begin{equation} \label{eq:Ns}
\sum_{\delta= 0}^{\infty}\overline{N}_{k,[S]}(\delta,LK_S)q^\delta = (\widetilde {DG}_2(y,q)/q)^k 
\frac{ B_1(y,q)^{K_S^2}B_2(y,q)^{LK_S} (D \widetilde {DG}_2(y,q)/q)}{(\widetilde \Delta(y,q)\cdot D\widetilde {DG}_2(y,q)/q^2)^{\chi(\oo_S)/2}}.
\end{equation}

\vspace{2mm} 
When $S$ is a K3 surface, the RHS simplifies dramatically.  Moreover, on the LHS, for each term
we may choose a representative K3 surface $S_g$ of genus $g$, with an irreducible line bundle $L_g$ giving the polarization.  That is,
$\overline{M}_{k,[K3]}(g-1, 0) = \overline{N}^{g-k}_{g-k,[S_g,L_g]}$.  The relative Hilbert schemes 
over the general $\proj^\delta \subset |L_g|$ are all smooth \cite{Muk}, so the 
$\overline{N}^\delta_{\delta,[S_g,L_g]}$ are equal to the geometric $\overline{N}^{\delta}_{\cC/\proj^\delta}$. 
Summarizing the preceding discussion,

\begin{conj} \label{conj:GconjK3}  For any $k$, 
\[ \sum_{g=k}^\infty q^{g-1} \overline{N}^{g-k}_{g-k,[S_g,L_g]}(y)  = \frac{\widetilde{DG}_2(y,q)^k}{\widetilde{\Delta}(y,q)}.\]
\end{conj} 

\begin{prop}
	Assume Conjecture \ref{DLconj}.  Then Conjectures \ref{Gconj} and \ref{conj:GconjK3} are equivalent.
\end{prop}
\begin{proof}
	According to Proposition \ref{produc}, from Conjecture \ref{DLconj} we can deduce that the 
	$\overline{N}$ have a multiplicative generating series.  The two of its components which are made
	explicit in Conjecture \ref{Gconj} are determined by the case when $S$ is a K3 surface, and we have
	seen that the posited explicit formula is equivalent to that given in Conjecture \ref{conj:GconjK3}. 
\end{proof}

\begin{prop}
Conjecture \ref{conj:GconjK3} is true at $k = 0$.
\end{prop}
\begin{proof}
	The relative Hilbert schemes over the general $\proj^\delta \subset |L_g|$ are smooth, 
	so the quantity $\overline{N}^{g}_{g,[S_g,L_g]}$ is the $\chi_{-y}$ genus of the relative compactified
	Jacobian of the tautological family of curves $\cC/|L_g|$.  In fact, Kawai and Yoshioka compute the 
	Hodge polynomial of this space \cite{KY}.  Alternatively, as in \cite{Bea} one may note that the relative
	compactified Jacobian over $|L_g|$ is birational to the Hilbert scheme of points $S_g^{[g]}$; since 
	both are hyperK\"ahler they are deformation equivalent by \cite{Hu}, and the Hodge polynomial
	of the Hilbert scheme of points was computed in \cite{GS}.  In either case the result is 
	given by the RHS. 
\end{proof}

We do not know how to compute the $\chi_{-y}$ genera of relative Hilbert schemes over linear subsystems for
K3 surfaces.  In fact, while the Euler numbers of these spaces are known, the only known calculation of them
(due to Maulik, Pandharipande, and Thomas) is 
extremely indirect, and involves at least two uses of the Gromov-Witten/Pairs correspondence \cite{MPT}.  
Rather remarkably, the same series which we have conjectured describes the 
$\overline{N}^\delta_\delta$ is found in \cite{MPT} to encode {\em all} the Euler numbers
of relative Hilbert schemes over linear subsystems. 

\begin{thm}\cite{MPT} \label{thm:MPT}  
For each $g$, let $(S_g,L_g)$ be a K3 surface of genus $g$, and assume $L_g$ is irreducible.
Let $H$ be the hyperplane class on $|L_g|$.  Then for any $k$, 
\[
(y-2+y^{-1})^{k-1} \frac{\widetilde{DG}_2^k}{\widetilde{\Delta}} =
\sum_{g=k}^\infty q^{g-1} \sum_{n=0}^\infty y^{n+1-g} \int_{\cC^{[n]}_{|L_g|}} c_{n+g-k}(T \cC^{[n]}_{|L_g|})
	\cdot \rho^*(H^k). \\
\]
\end{thm}

Comparing powers of $q$, we see that given \cite{MPT}, Conjecture \ref{conj:GconjK3} is equivalent
to the following statement: 

\begin{conj} \label{conj:tryresidue}
  Let $(S, L)$ be a $K3$ surface of genus $g$ with $L$ irreducible.  For all $k$, 
  \[ 
  (y - 2 + y^{-1})^{k-1} \overline{N}^{g-k}_{g-k, [S,L]} = \sum_{n=0}^\infty 
  y^{n+1-g} \int_{\cC^{[n]}_{|L|}} c_{n+g-k}(T \cC^{[n]}_{|L|})\cdot \rho^*(H^k).
  \]
\end{conj}

\begin{rem}
The statement of Theorem \ref{thm:MPT} 
in \cite{MPT} differs slightly; there are some signs owing to the use of $\Omega$ rather than
$T$, and it is formulated in terms of the space of ``stable pairs'' rather than the relative Hilbert scheme. 
But since $L_g$ is irreducible, these are the same space as per \cite[Appendix B]{PT2}.  In the language
of \cite{MPT}, Conjecture
\ref{conj:tryresidue} may be viewed as asserting that the $\overline{N}$ encode certain descendent
integrals in the stable pairs theory, or equivalently in the Gromov-Witten theory.  
\end{rem}

We may instead specialize to abelian surfaces.  Let $(A, L_g)$ denote a primitively polarized abelian surface
with $L_g^2 = 2g-2$ and hence $\chi(L) = g-1$ and $g(L) = g$; assume $L$ is irreducible.  Note such surfaces
exist for all $g \ge 1$.  
Equation \ref{eq:Ms} specializes to

\[\sum_{g=k+2}^\infty \overline{N}^{g-k-2}_{g-k-2,[A,L_g]} q^{g-1}  = \sum_{l=k+1}^\infty \overline{M}_{k,[A]}(l,0)q^l  = \widetilde{DG}_2(y,q)^k D\widetilde{DG}_2
(y,q).\] 

The above formula is not equivalent to Conjecture \ref{Gconj}, since $\chi(\mathcal{O}_A) = 0$.  Moreover we do
not know how to establish it, even at $k = 0$.  However:

\begin{prop}
Assume Conjecture \ref{DLconj} or $K_S = 0$. Conjecture \ref{Gconj} is equivalent to the following two formulas: 
\begin{eqnarray}
\label{GformK3} \sum_{g=0}^\infty \overline{N}^g_{g,[K3,L_g]} q^{g-1} & = & \widetilde{\Delta}(y,q)^{-1}, \\
\label{GformA} \sum_{g=2}^\infty \overline{N}^{g-2}_{g-2,[A,L_g]} q^{g-1} & = & D \widetilde{DG}_2(y,q).
\end{eqnarray}
\end{prop}
\begin{proof}
  It is enough to show that the given invariants of the complete linear system suffice to determine 
  the series $A_1, A_4$ and then apply the residue trick explained above.  But
  $\overline{N}^{g-2}_{g-2,[A,L_g]}(y) = \Coeff_{s^{g-2}} A_1(y,s)^{2g-2}$.  Beginning at $g = 2$ this allows
  to iteratively determine the coefficients of $A_1$.  Then 
  $\overline{N}^g_{g,[K3,L_g]}(y) = \Coeff_{s^g} A_1(y,s)^{2g-2} A_4(y,s)^{24}$; since we now know $A_1$ this
  determines $A_4$.  
\end{proof}

As we have established Equation \ref{GformK3}, it remains only to compute the invariants for complete linear
systems on abelian surfaces. In the meantime in \cite{GSh} we have proved a generalized version of Equation \ref{GformA} using the methods of \cite{KY}.
By the above this gives a proof of 
\conjref{Gconj} in the case $K_S$ numerically trivial. 
This also provides an entirely sheaf theoretic proof of the formula \cite[Thm. 6]{MPT} for linear subsystems for irreducible
line bundles on K3 surfaces.

\section{Refined Severi degrees} \label{sec:Severi}

The {\em Severi degrees} $n^{d,\delta}$ are the numbers of $\delta$-nodal reduced 
degree $d$ curves in $\P^2$ though $\binom{d+2}{2}-1-\delta$ general points. The 
famous Caporaso-Harris formula \cite{CH} gives a recursive method of computing the 
Severi degrees. The recursion 
involves the {\em relative Severi degrees} $n^{d,\delta}(\alpha,\beta)$ (we sketch their definition below) which count 
delta-nodal curves with tangency conditions along a fixed line in $\P^2$.
More generally, for a line bundle $L$ on a surface $S$, one can define the {\em Severi degree}
$n^{L,\delta}$ as the number of $\delta$-nodal reduced curves in $|L|$ though $\dim |L|-\delta$ general 
points, provided this number is finite.

We begin by a review of the Caporaso-Harris recursion formula; we will use the more general formulation of Vakil \cite{V}, which also applies to rational ruled surfaces. 
By a {\it sequence} we mean a collection $\alpha=(\alpha_1,\alpha_2,\ldots)$ of nonnegative integers, almost all of which are zero. We write $d$ for the sequence $(d,0,0,\ldots)$ and $e_k$ for the sequence whose $k$-th element is $1$ and all other ones  $0$.
For two sequences
$\alpha$, $\beta$ we define 
$|\alpha|=\sum_i\alpha_i$, $I\alpha=\sum_i i\alpha_i$, 
$\alpha+\beta=(\alpha_1+\beta_1,\alpha_2+\beta_2,\ldots)$, $\binom{\alpha}{\beta}=\prod_i \binom{\alpha_i}{\beta_i}$. We write $\alpha\le\beta$ to mean $\alpha_i\le \beta_i$ for all $i$.

Throughout this section we take $S$ to be $\P^2$ or a rational ruled surface. 
In case $S=\P^2$, let $E$ be a line in $\P^2$, in case $S$ is a rational ruled surface $\Sigma_e=\P(\oo_{\P^1}\oplus\oo_{\P^1}(-e))$, let $E$ be the class of the section with $E^2=-e$.
We denote $H$ the hyperplane class on $\P^2$, $F$ the class of a fibre on $\Sigma_e$.

Let  $L$ be a line bundle on $S$ and let $\alpha$, $\beta$ be sequences with $I\alpha+I\beta=EL$, and 
let $\delta\ge 0$ be an integer. Let $\gamma(L, \beta,\delta):=\dim|L|-EL+|\beta|-\delta$.
 The  {\it relative Severi degree} $n^{L,\delta}(\alpha,\beta)$ is the number of $\delta$-nodal curves $C$ in $|L|$, which do not contain $E$ as a component and for each $k$ with $\alpha_k$  points of contact of order $k$ to $E$ at given points of $E$ and $\beta_k$ points of contact of order $k$, to $E$ at variable points of $E$, and passing through 
 $\gamma(L, \beta,\delta)$ general points of $S$ (see \cite{V} for a more formal definition).

\begin{rec}\label{CHrec} \cite{CH,V}
The  relative Severi degrees 
$n^{L,\delta}(\alpha,\beta)$ are 
recursively given  as follows: 
$n^{L,\delta}(\alpha,\beta)=0$ if $\gamma(L, \beta,\delta)<0$.
If $\gamma(L, \beta,\delta)>0$, then
\begin{equation}\label{Srec}
\begin{split}
n^{L,\delta}&(\alpha,\beta)=\sum_{k:\beta_k>0} k \cdot n^{L,\delta}(\alpha+e_k,\beta-e_k)
\\&+\sum_{\alpha',\beta',\delta'}
\prod_i i^{\beta_i'-\beta_i}\binom{\alpha}{\alpha'}\binom{\beta'}{\beta} n^{L-E,\delta'}(\alpha',\beta').
\end{split}
\end{equation}
Here the second sum runs through all  $\alpha',\beta',\delta'$ satisfying the condition
\begin{equation}\label{relcong}
\begin{split}
\alpha'&\le \alpha, \ \beta'\ge \beta,\ I\alpha'+I\beta'=E(L-E),\\ \delta'&=\delta+g(L-E)-g(L)+|\beta'-\beta|+1=\delta-E(L-E)+|\beta'-\beta|.
\end{split}
\end{equation}
{\bf Initial conditions:} if $\gamma(L,\beta,\delta)=0$ we have  $n^{L,\delta}(\alpha,\beta)=0$ unless we are in one of the following cases
\begin{enumerate}
\item 
In case $S=\P^2$ we put $n^{H,0}(1,0)=1$,
\item
In case $S=\Sigma_e$, let $F$ be the class of a fibre of the ruling; we put $n^{kF,0}(k,0)=1$.
\end{enumerate}
\end{rec}
We put $n^{L,\delta}:=n^{L,\delta}(0,LE)$. In case $S=\P^2$, we write $n^{d,\delta}(\alpha,\beta):=n^{dH,\delta}(\alpha, \beta),$ and $n^{d,\delta}:=n^{dH,\delta}(0, d),$

\subsection{Refined Severi degrees}

We formally introduce a refinement of this recursion. 

\begin{rec} \label{refCHrec}
With the same notations, assumptions, limits of summation, and initial values as for the relative Severi degrees in 
Recursion \ref{CHrec}, we define the refined relative Severi degrees $N^{L,\delta}(\alpha,\beta)(y)$ for 
$\gamma(L, \beta,\delta)>0$ by:
\begin{equation}\label{refrec}
\begin{split}
N^{L,\delta}&(\alpha,\beta)(y)=\sum_{k:\beta_k>0} \frac{1-y^k}{1-y} \cdot N^{L,\delta}(\alpha+e_k,\beta-e_k)(y)
\\&+\sum_{\alpha',\beta',\delta'}
y^{I\alpha'+I\beta}\prod_i \left(\frac{1-y^i}{1-y}\right)^{\beta_i'-\beta_i}\binom{\alpha}{\alpha'}\binom{\beta'}{\beta} N^{L-E,\delta'}(\alpha',\beta')(y).
\end{split}
\end{equation}
\end{rec}

%

We abbreviate $N^{L,\delta}:=N^{L,\delta}(0,LE)$,  and,  
in case $S=\P^2$, $N^{d,\delta}(\alpha, \beta):=N^{dH,\delta}(\alpha, \beta)$, 
$N^{d,\delta}:=N^{dH,\delta}(0,d)$.
As with the refined invariants, we define normalized refined relative Severi degrees  which are 
Laurent polynomials in $y^{1/2}$, symmetric under $y\mapsto 1/y$.

\begin{defn}\label{barNL}
The {\it normalized (relative) refined Severi degrees} 
$\overline N^{L,\delta}(y)$ are defined by 
$$\overline N^{L,\delta}(\alpha,\beta)(y)=N^{L,\delta}(\alpha,\beta)(y)/y^{\delta+(I\beta-|\beta|)/2}, \quad
\overline N^{L,\delta}(y)=N^{L,\delta}/y^{\delta}.$$
\end{defn}

\begin{prop}
The $\overline N^{L,\delta}(\alpha,\beta)(y)$ are determined by the same initial conditions as the 
$N^{L,\delta}(\alpha,\beta)(y)$ and the recursion
\begin{equation}\label{norrefrec}
\overline N^{L,\delta}(\alpha,\beta)=\sum_{k:\beta_k>0} [k]_y \cdot \overline N^{L,\delta}(\alpha+e_k,\beta-e_k)
+\sum_{\alpha',\beta',\delta'}
\left(\prod_i [i]_y^{\beta_i'-\beta_i}\right)
\binom{\alpha}{\alpha'}\binom{\beta'}{\beta} \overline N^{L-E,\delta'}(\alpha',\beta')
\end{equation}
with the same conditions on $\alpha',\beta',\delta'$ as above.
In particular $\overline N^{L,\delta}(\alpha,\beta)(y)$ is symmetric under $y\mapsto 1/y$.
\end{prop}
\begin{proof}
It is enough to prove that every summand on the right hand side of \eqref{refrec} is obtained from the corresponding summand of \eqref{norrefrec} by multiplying by $y^{\delta+(I\beta-|\beta|)/2}$.
Each summand in the first sum is multiplied by $y^m$ with $$m=(k-1+I(\beta-e_k)-(|\beta-e_k|))/2+\delta=\delta+(I\beta-|\beta|)/2.$$
Each summand in the second sum is multiplied by $y^m$ with 
$$m=I\alpha'+I\beta+(I(\beta'-\beta)-|\beta'-\beta|)/2+\delta'+(I\beta'-|\beta'|)/2=\delta+(I\beta-|\beta|)/2,$$
where we use $I\alpha'+I\beta'=E(L-E)=\delta-\delta'+|\beta'-\beta|$.
\end{proof}

It is clear that the recursions for the refined Severi degrees specialize at $y=1$ to the recursion for the usual
Severi degrees.  Thus,

\begin{prop} 
$N^{L,\delta}(\alpha,\beta)(1)=\overline N^{L,\delta}(\alpha,\beta)(1)=n^{L,\delta}(\alpha,\beta)$.
\end{prop}

According to \cite{KS}, if the general $\proj^\delta \subset |L|$ contains no non-reduced curves
and no curves containing components with negative self intersection, the Severi degrees are computed
by the universal formulas: $n^{L,\delta} = n^\delta_{\delta,[S,L]}$.  We expect the same for
refined Severi degrees. 

\begin{conj}\label{ref-sev}
Let $S$ be $\P^2$ or a rational ruled surface, let $L$ be a line bundle, and assume 
$\proj^\delta \subset |L|$ contains no non-reduced curves
and no curves containing components with negative self intersection.  
Then the refined Severi degrees are computed
by the universal formulas: $N^{L,\delta} = N^\delta_{\delta,[S,L]}$.  Explicitly,
\begin{enumerate}
\item On $\P^2$ we have $N^{d,\delta}=N^{\delta}_{\delta,[\P^2,dH]}$, for $\delta\le 2d-2$.
\item On $\P^1\times\P^1$ we have $N^{nF+mG,\delta}=N_{\delta,[\P^1\times \P^1,nF+mG]}^\delta$ 
for $\delta\le 
\min(2n,2m)$. 
\item On $\Sigma_e$ with $e>0$, we have $N^{nF+mE,\delta}=N^\delta_{\delta,[\Sigma_e,nF+mE]}$ 
for $\delta\le \min(2n,n-em)$.
\end{enumerate}
\end{conj}

Directly from the defining Recursion \ref{refCHrec} we have computed all the  $N^{d,\delta}(y)$ for $d\le 15$ and $\delta\le 30$. Assuming the vanishing Conjecture \ref{DLconj} and part (1) of \conjref{ref-sev}, the refined 
Severi degrees suffice to determine
all the power series in Corollary \ref{produc} 
or equivalently in Conjecture \ref{Gconj}.  Note that Recursion \ref{refCHrec}
is much more computationally tractable than equivariant localization.  Thus under the above assumption, 
we have verified \conjref{Gconj} modulo $q^{29}$ and determined $B_1(y,q)$ and $B_2(y,q)$ modulo $q^{29}$.

\subsection{Irreducible refined Severi degrees}
Denote by $n_0^{L,\delta}$ the irreducible Severi degrees, i.e. informally the number of irreducible $\delta$-nodal curves in $|L|\ne |E|$ passing though $\dim|L|-\delta$ general points. 
In \cite{Ge} Getzler observes in the case $S=\P^2$  that the $n_0^{d,\delta}$ can be expressed in terms of the 
Severi degrees $n^{d,\delta}$ by the relation
\begin{align*}
\sum_{d,\delta} \frac{z^{\binom{d+2}{2}-\delta-1}}{(\binom{d+2}{2}-\delta-1)!} q^d n_0^{d,\delta}=
\log\left(1+\sum_{d,\delta} \frac{z^{\binom{d+2}{2}-\delta-1}}{(\binom{d+2}{2}-\delta-1)!} q^d n^{d,\delta}\right)\!.
\end{align*}
The generalization of this to $n_0^{L,\delta}$ is in \cite{V}. The same formula can be used to define the  irreducible normalized refined Severi degrees $\overline N_0^{L,\delta}(y)$ by 
\begin{align*}
\sum_{L,\delta} \frac{z^{\dim|L|-\delta}}{(\dim|L|-\delta)!} v^L \overline N_0^{L,\delta}(y)=
\log\left(1+\sum_{L,\delta}\frac{z^{\dim|L|-\delta}}{(\dim|L|-\delta)!} v^L  \overline N^{L,\delta}(y)\right)\!,
\end{align*}
and the irreducible refined Severi degrees by $N_0^{L,\delta}(y):=y^\delta\overline N_0^{L,\delta}(y)$.
Here $\big\{v^L\big\}_{L \text{effective}, L\ne E}$ are elements of the Novikov ring, i.e. $v^{L_1}v^{L_2}=v^{L_1+L_2}$.
Evidently $\overline N_0^{L,\delta}(y)$ is a Laurent polynomial in $y$ invariant under $y\mapsto 1/y$, and
$N_0^{L,\delta}(1)=\overline N_0^{L,\delta}(1)=n_0^{L,\delta}$.

\begin{thm} \cite{BG}
$N_0^{L,\delta}$ has nonnegative integer coefficients. 
\end{thm}

From this positivity, one can conclude vanishing results for $N_0^{L,\delta}$ from the analogous (known)
results for $n_0^{L,\delta}$.  For instance $N_0^{L,\delta}(y)=n_0^{L,\delta}=0$ for $\delta>g(L)$
since there are no irreducible curves of cogenus greater than $g(L)$. 

We list the first few of the $N_0^{d,\delta}(y)$.
Write $N_0^d(y,t):=\sum_{\delta\ge 0}N_0^{d,\delta}(y)t^\delta$.  We have computed the  $N_0^d(y,t)$  for 
$d\le 14$.
We have
\begin{align*}N_0^1(y,t)&=1,\quad N_0^2(y,t)=1, \quad, N_0^3(y,t)=1+(y^2 + 10y + 1)t, \\
N_0^4(y,t)&=1 + (3y^2 + 21y + 3)t + (3y^4 + 33y^3 + 153y^2 + 33y + 3)t^2\\&  + (y^6 + 13y^5 + 94y^4 + 404y^3 + 94y^2 + 13y + 1)t^3,\\
N_0^5(y,t)&=1 + (6y^2 + 36y + 6)t + (15y^4 + 156y^3 + 540y^2 + 156y + 15)t^2 \\
&+ (20y^6 + 268y^5 + 1555y^4 + 4229y^3 + 1555y^2 + 268y + 20)t^3 \\
&+ (15y^8 + 228y^7 + 1674y^6 + 7407y^5 + 18207y^4 + 7407y^3 \\&
+ 1674y^2 + 228y + 15)t^4 + (6y^{10} + 96y^9 + 792y^8 + 4398y^7\\& + 17190y^6 + 42228y^5 + 17190y^4 + 4398y^3 + 792y^2 + 96y + 6)t^5\\& + (y^{12} + 16y^{11} + 139y^{10} + 867y^9 + 4203y^8 + 16377y^7 \\&+ 44098y^6 + 16377y^5 + 4203y^4 + 867y^3 + 139y^2 + 16y + 1)t^6.
\end{align*}
On $\P^1\times\P^1$ we have computed the $N_0^{nF+mG,\delta}$ for $1\le, n,m\le 8$. 
We list the first few; write $N_0^{n,m}:=\sum_{\delta\ge 0}N_0^{n,m,\delta}(y)t^\delta$. We have
\begin{align*}
N^{1,k}_0(y,t)&=1, \hbox{ for all $k$}, \
N^{2,2}_0(y,t) = 1 + (y^2 + 10y + 1)t ,\\
N^{2,3}_0(y,t) &= 1 + (2y^2 + 16y + 2)t + (y^4 + 12y^3 + 79y^2 + 12y + 1)t^2,\\
N^{2,4}_0 (y,t)&= 1 + (3y^2 + 22y + 3)t + (3y^4 + 36y^3 + 174y^2 + 36y + 3)t^2 \\
&+ (y^6 + 14y^5 + 117y^4 + 596y^3 + 117y^2 + 14y + 1)t^3,\\
N^{3,3}_0(y,t) &= 1 + (4y^2 + 26y + 4)t + (6y^4 + 64y^3 + 256y^2 + 64y + 6)t^2\\
& + (4y^6 + 52y^5 + 332y^4 + 1168y^3 + 332y^2 + 52y + 4)t^3 \\
+ (y^8 + 14y^7 &+ 109y^6 + 636y^5 + 2430y^4 + 636y^3 + 109y^2 + 14y + 1)t^4.
\end{align*}

\subsection{The refined invariants at $y=0$}

Now we compute the specialization of the refined Severi degrees at $y=0$.

\begin{nota} For a sequence $\beta$ we write $\binom{|\beta|}{\beta}:=\frac{|\beta|!}{\prod_i \beta_i!}.$
\end{nota}

\begin{prop}\label{chi0}
$N^{L,\delta}(\alpha,\beta)(0)=\binom{|\beta|}{\beta}\binom{g(L)}{\delta}.$
In particular $N^{L,\delta}(0)=\binom{g(L)}{\delta}$.
\end{prop}
\begin{proof}
It is easy to see that the statement holds for the initial values.
Setting $y=0$ in the recursion formula \eqref{refrec} gives the two recursion formulas
\begin{equation}
\label{rec1}
\begin{split}
N^{L,\delta}(\alpha,\beta)(0)&=\sum_{k|\beta_k>0} N^{L,\delta}(\alpha+e_k,\beta-e_k)(0), \quad\hbox{if $\beta\ne 0$},\\
N^{L,\delta}(\alpha,0)(0)&=\!\!\!\!\! \sum_{I\beta=E(L-E)}  \!\!\!\!\! N^{L-E,\delta-E(L-E)+|\beta|}(0,\beta)(0).
\end{split}
\end{equation}
The first formula gives 
\begin{equation}\label{rec3}
N^{L,\delta}(\alpha, \beta)(0)=\binom{|\beta|}{\beta}N^{L,\delta}(\alpha+\beta,0)(0).
\end{equation}
The second formula shows in particular that $N^{L,\delta}(\alpha,0)(0)$ is independent of $\alpha$, thus 
$N^{L,\delta}(\alpha,0)(0)=N^{L,\delta}(LE,0)(0)=N^{L,\delta}(0)$; the last equality is by \eqref{rec3}.
Thus $N^{L,\delta}(\alpha, \beta)(0)=\binom{|\beta|}{\beta}N^{L,\delta}(0).$

Finally, writing $L=L_0+aE$ with $L_0=H$ in case $S=\P^2$ and $L_0=kF$ in case $S$ is a rational ruled surface, we prove $N^{L,\delta}=\binom{g(L)}{\delta}$ by induction over $a$. It is easy to see that the claim is true for $L_0$. 
By induction the second equation of \eqref{rec1} becomes
$$N^{L,\delta}(0)=\sum_{I\beta=E(L-E)} \binom{g(L-E)}{\delta-E(L-E)+|\beta|}\binom{|\beta|}{\beta}.$$
Thus we need to show the identity
\begin{equation}\label{rec4} (1+t)^{g(L)}=\sum_{\delta\ge 0} \sum_{I\beta=E(L-E)} \binom{g(L-E)}{\delta-I\beta+|\beta|}\binom{|\beta|}{\beta}t^\delta.
\end{equation}
Note that by the  multinomial formula we have
\begin{align*}\frac{x}{1-x(1+t)}&=\frac{1}{1-\frac{x}{1-tx}}-1=\sum_{n>0} (x+tx^2+t^2x^3+\ldots)^n\\&=\sum_{\beta\ne 0}
\binom{|\beta|}{\beta} t^{I\beta-|\beta|}x^{I\beta}.
\end{align*}
Thus the right hand side of \eqref{rec4} becomes
$$\Coeff_{x^{E(L-E)}}\left[(1+t)^{g(L-E)} \frac{x}{1-x(1+t)}\right]=(1+t)^{g(L)}.$$
\end{proof}

Using the tropical interpretation of the $N^{L,\delta}(\alpha,\beta)$ of \cite{BG} using refined multiplicities (see also Section \ref{sec:RRT}), this result has been generalized in \cite{IM} to arbitrary toric surfaces.

For $M$ a line bundle on $S$, let $M_n:=f^*g_*(\otimes_{i=1}^n pr_i^*M)^{{\mathfrak S}_n}\in \mathrm{Pic}\,
S^{[n]}$, 
where 
$f:S^{[n]}\to S^{(n)}$ and $g:S^n\to S^{(n)}$ are the natural morphisms, and $pr_i:S^n\to S$ is the $i^{th}$ projection. It is well-known that $K_{S^{[n]}}=(K_S)_n$, and it is also standard that $\det M^{[n]}=M_n\otimes \det \oo_S^{[n]}$ (see e.g. \cite{EGL}).
\begin{lem}\label{Sca}
$$\sum_{n,k\ge 0} \chi(S^{[n]},\Lambda^k(L^{[n]})^\vee)x^kq^n=
\frac{(1+xq)^{\chi(L^\vee)}}{(1-q)^{\chi(\oo_S)}}.$$
\end{lem}
\begin{proof}
This is a corollary to \cite[Thm.~5.2.1]{Sc}, which implies for line bundles $L,M$ on $S$ that
$\chi(S^{[n]},\Lambda^kL^{[n]}\otimes M_n)=\binom{\chi(L\otimes M)}{k}\binom{\chi(M)+n-k-1}{n-k}$.
We apply this with $M=K_S$.
Thus we have by applying Serre duality on $S^{[n]}$ and on $S$:
\begin{align*}\chi(&S^{[n]},\Lambda^k(L^{[n]})^\vee)=\chi(S^{[n]},\Lambda^kL^{[n]}\otimes (K_S)_n)\\&=
\binom{\chi(L\otimes K_S)}{k}\binom{\chi(K_S)+n-k-1}{n-k}=\binom{\chi(L^\vee)}{k}\binom{\chi(\oo_S)+n-k-1}{n-k},\end{align*}
which is equivalent to the statement of the Lemma.
\end{proof}

\begin{prop}\label{chi01}
$N^l_{\delta,[S,L]}(0)=0$ for $l>\delta$, and for all $0\le l\le \delta$ we have
$N^l_{\delta,[S,L]}(0)=\binom{\chi(L^\vee)}{l}=\binom{g(L)+\chi(\oo_S)-1}{l}$.
In particular if $S$ is a rational surface, then $N^\delta_{\delta,[S,L]}(0)=\binom{g(L)}{\delta}$.
\end{prop}
\begin{proof}
By  \propref{prop:hilbgenus} and \defref{Ndef}, we have
$$\Coeff_{q^n}\left[\sum_{l\ge 0} N^{l}_{\delta,[S,L]}(0)\frac{q^l}{(1-q)^{l+1-g(L)}}\right] =\res_{x=0}\int_{S^{[n]}}
\left(\frac{1}{1-e^{-x}}\right)^{\delta+1}\frac{c_n(L^{[n]}\cdot e^x)td(S^{[n]})}{td(L^{[n]}\cdot e^x)}dx.$$
Note that by definition $\frac{c_n(L^{[n]}\cdot e^x)}{td(L^{[n]}\cdot e^x)}=\sum_{k=0}^n (-e^{-x})^k
\ch(\Lambda^k(L^{[n]})^\vee).$
Thus by  Riemann-Roch the right hand side is
$$\res_{x=0}\left[\left(\frac{1}{1-e^{-x}}\right)^{\delta+1}\sum_{k=0}^n (-e^{-x})^k
\chi(S^{[n]},\Lambda^k(L^{[n]})^\vee)dx\right].$$
We put $T=e^{-x}$ and apply \lemref{Sca} to obtain
\begin{align*}
\sum_{l\ge 0} N^{l}_{\delta,[S,L]}(0)\frac{q^l}{(1-q)^{l+1-g(L)}}&=-\res_{T=1}\left[\left(\frac{1}{1-T}\right)^{\delta+1}\sum_{n\ge 0}\sum_{k=0}^n (-T)^kq^n
\chi(S^{[n]},\Lambda^k(L^{[n]})^\vee)\frac{dT}{T}\right]\\
&=-\res_{T=1}\left[\left(\frac{1}{1-T}\right)^{\delta+1}\frac{(1-Tq)^{\chi(L^\vee)}}{(1-q)^{\chi(\oo_S)}} 
\frac{dT}{T}\right].\end{align*}
Substituting 
$T=1-\alpha$ the right hand side becomes
\begin{align*}&\res_{\alpha=0}\left[\frac{1}{\alpha^{\delta+1}}
\frac{(1-q+\alpha q)^{\chi(L^\vee)}}{(1-q)^{\chi(\oo_S)}}\sum_{l\ge 0}\alpha^l d\alpha 
\right]=\sum_{l=0}^\delta \binom{\chi(L^\vee)}{l} q^l(1-q)^{\chi(L^\vee)-l-\chi(\oo_S)}\\&=
\sum_{l=0}^\delta \binom{\chi(L^\vee)}{l} q^l(1-q)^{g(L)-l-1}.
\end{align*}
\end{proof}

\subsection{Conjectural generalization to higher powers of $y$}
\propref{chi01} and \propref{chi0} can be subsummed in the following statements:
\begin{enumerate}
\item For any line bundle $L$ on a surface $S$ we have 
$$\sum_{\delta\ge 0} N^\delta_{\delta,[S,L]}(0) q^\delta=(1+q)^{\chi(L^\vee)}.$$
\item If  $L$ is an effective  divisor on $\P^2$ or a rational ruled surface $S$, then 
 $$\sum_{\delta\ge 0} N^\delta_{\delta,[S,L]}(0)N^{L,\delta}(0) q^\delta=(1+q)^{g(L)}=\sum_{\delta\ge 0} N_{\delta,[S,L]}^\delta(0) q^\delta.$$
\end{enumerate}
We want to give a conjectural extension of these two statements to higher powers of $y$.
We start with the analogue of (1)
\begin{conj}\label{yconj}
Let $L$ be a line bundle on a surface $S$.
Then we have for all $i\ge 0$
$$\Coeff_{y^i}\left[\sum_{\delta\ge 0} N^\delta_{\delta,[S,L]}(y) q^\delta\right]
=(1+q)^{\chi(L^\vee)-3i}P^i_{L}(q).$$
Here $P^i_{L}(q)$  is a polynomial in $q$ of degree at most $3i$.
In particular if $\chi(L^\vee)\ge 3i$ then $\Coeff_{y^i}N^\delta_{\delta,[S,L]}(y)=0$ for $\delta>\chi(L^\vee)$. 

\end{conj}

Assuming \conjref{DLconj} we get by \remref{produc} 
$$\sum_{\delta\ge 0} N_{\delta,[S,L]}^\delta(y) q^\delta=F_1^{\chi(L^\vee)}F_2^{LK_S/2}F_3^{K_S^2}F_4^{\chi(\oo_S)},
$$
 with $F_i\in \Q[y]\lbr q\rbr $
We put
$C_1=\big(F_1/(1+q)\big)\big|_{y\mapsto y(1+q)^3}$, $C_2=F_2\big|_{y\mapsto y(1+q)^3}$, $C_3=F_3\big|_{y\mapsto y(1+q)^3}$, $C_4=F_4\big|_{y\mapsto y(1+q)^3}$.
\begin{conj}\label{ciconj}
For $i=1,\ldots, 4$ we have $C_i\in \Q[y^{1/3}]\lbr y^{1/3}q\rbr \cap \Q[q]\lbr y\rbr $.
\end{conj}

\begin{prop}
Assuming \conjref{ciconj},
 \conjref{yconj} holds. Furthermore
$$\sum_{i\ge 0}P^i_{L}(q) y^i=C_1^{\chi(L^\vee)}C_2^{LK_S/2}C_3^{K_S^2}C_4^{\chi(\oo_S)}.$$
\end{prop}
\begin{proof}
By definition 
$$(1+q)^{\chi(L^\vee)}C_1^{\chi(L^\vee)}C_2^{LK_S}C_3^{K_S^2}C_4^{\chi(\oo_S)}=
\sum_{\delta\ge 0} N_{\delta,[S,L]}^\delta(y(1+q)^3) q^{\delta}.$$
Therefore
$$\Coeff_{y^i}\left[ \sum_{\delta\ge 0} N_{\delta,[S,L]}^\delta(y) q^{\delta}\right]=(1+q)^{\chi(L^\vee)-3i}\Coeff_{y^i}\left[C_1^{\chi(L^\vee)}
C_2^{LK_S/2}C_3^{K_S^2}C_4^{\chi(\oo_S)}\right].$$ As by \conjref{ciconj}
all $C_i$ are in $\Q[y^{1/3}]\lbr y^{1/3}q\rbr $, we see that the coefficient of $y^i$ is a polynomial of degree at most $3i$ in $q$.
\end{proof}

\conjref{ciconj} has been verified modulo $q^{11}$. 
Assuming  \conjref{ref-sev} it has been verified modulo $q^{29}$.
We list the power series $C_1$, $C_2$, $C_3$, $C_4$ modulo $y^4$. 
\begin{align*}
C_1&=1+ (4q+2q^2)y + (q - 7q^2 + 12q^3 + 15q^4 + 3q^5)y^2 \\& + (-6q^2 + 56q^3 - 104q^4 - 112q^5 + 26q^6 + 32q^7 + 4q^8)y^3 + O(y^4),\\
C_2&=1 + (-2q - 6q^2 - 2q^3)y + (5q^2 + 48q^3 + 35q^4 + 6q^5 + q^6)y^2 \\&+ (14q^3 - 390q^4 - 286q^5 + 60q^6 + 52q^7)y^3 + O(y^4),
\\
C_3&=1 + (-q - 3q^2 - q^3)y + (q^2 + 16q^3 + 2q^4 - 6q^5 - q^6)y^2 \\&+ (15q^3 - 130q^4 + 66q^5 + 199q^6 + 65q^7)y^3 + O(y^4),\\
C_4&=1 + (6q + 18q^2 + 10q^3)y + (18q^2 + 64q^3 + 219q^4 + 222q^5 + 67q^6)y^2 \\ &+ (-44q^3 + 336q^4 + 72q^5 + 952q^6 + 2328q^7 + 1608q^8 + 352q^9)y^3 + O(y^4).
\\
\end{align*}

Now we formulate the conjectural analogue of (2). For simplicity we only deal with the case of $\P^2$ and $\P^1\times \P^1$.
\begin{conj}
\begin{enumerate}
\item Let $S=\P^2$ and 
assume $d\ge i+2$, then 
$$\Coeff_{y^i}\left[\sum_{\delta\ge 0} N^{d,\delta}(y)q^\delta\right]=(1+q)^{g(dH)-3i}P^i_{dH}(q).$$
\item Let $S=\P^1\times \P^1$ and assume $n,m\ge i+1$, then 
$$\Coeff_{y^i}\left[\sum_{\delta\ge 0} N^{nF+mG,\delta}(y)q^\delta\right]=(1+q)^{g(nF+mG)-3i}P^i_{nF+mG}(q).
$$
\end{enumerate}
\end{conj}
For $d\le 14$, and for $n,m\le 8$ this conjecture has been checked  modulo $q^{11}$ and, assuming \conjref{ref-sev}, modulo $q^{29}$.

\section{Refined, real, and tropical} \label{sec:RRT}

Mikhalkin \cite{Mi} has shown that the Severi degrees of projective toric surfaces can also be computed  
using tropical geometry: the Severi degrees  $n^{L,\delta}$ count -- with multiplicities -- simple tropical curves through
$\dim|L|- \delta$ points in $\R^2$ in tropical general position. 
Roughly speaking a simple  tropical curve $C$  is a trivalent graph $\Gamma$ immersed in $\R^2$ together with some extra data.  From this data, one assigns to each vertex $v$ of $\Gamma$ 
 a multiplicity $m(v)\in \Z_{\ge 0}$ and defines 
the {\it complex multiplicity $m(C)$} as the product of the $m(v)$ over the vertices of $\Gamma$. 
In \cite{GM} a proof of the Caporaso-Harris recursion formula is given via tropical geometry. 

The analogues of the Gromov-Witten invariants in real algebraic geometry are the Welschinger invariants  \cite{W}.  These were originally
defined to count real pseudoholomorphic curves in real symplectic manifolds.
We restrict attention to the case that $S$ is a smooth projective toric surface.  As toric varieties are defined
over $\mathbb{Z}$, they certainly carry a real structure, and we write $\sigma$ for the associated anti-holomorphic
involution.  A real curve in $S$ is  an algebraic curve $C\subset S$  with $C=\sigma(C)$, 
and the real locus of $C$ is $C^\sigma$.
Fix a generic set $\Sigma$ of $\dim |L|-\delta$ general real\footnote{
We are here only considering the so-called {\it totally real} Welschinger invariants. 
More generally one could consider for any $0\le l\le (\dim|L|-\delta)/2$ the numbers
$W^{L,\delta,l}(\Sigma)$ which count real curves passing through $\dim|L|-\delta-2l$ real  points and 
$l$ pairs of complex conjugate points. 
} points on $S$. The real enumerative invariant 
is $W^{L,\delta}(\Sigma):= \sum_{C} (-1)^{s(C)}$, where $C$ runs 
through the possibly reducible real curves $C\in |L|$ of geometric genus $g(L)-\delta$, 
passing through all the points of $\Sigma$, and $s(C)$ is the number of isolated real
 nodes of $C$, i.e. the points where $C$ analytically locally has the equation $x^2+y^2$. 
We denoted by $W_{0}^{L,\delta}(\Sigma)$ the corresponding sum for irreducible curves. 
If $S$ is an unnodal (i.e. it contains no rational curve with self intersection $-n$, with $n\ge 2$) del Pezzo surface then the real enumerative invariants coincide with the Welschinger invariants.
In \cite{W} it was proven that  $W_{0}^{L,g(L)}(\Sigma)$, i.e. the count of curves of geometric genus $0$, is independent of the generic $\Sigma$.
We will denote it just by $W_{0}^{L,g(L)}$. In general 
$W^{L,\delta}(\Sigma)$ and $W^{L,\delta}_0(\Sigma)$ will depend on $\Sigma$ via a system of walls and chambers.

In a sense, we have already seen these invariants.  For a family of real curves $\cC/B$, let
$n^{i, \R}_{\cC/B}$ be defined by the same formula as the $n^{i,\R}_C$ introduced for individual curves
in Section \ref{sec:curve}.  Then we have:

\begin{prop}
Let $L$ be a real line bundle on $S$, and let $\proj^\delta \subset |L|$ be a linear subsystem determined by the 
real point conditions $\Sigma$.  Assume that all curves in $\proj^\delta$ are reduced, that no curves have cogenus
greater than $\delta$, and that all curves of cogenus $\delta$ are nodal.  Then
\[(-1)^\delta n^{\delta, \R}_{\cC/\proj^\delta} = W^{L,\delta}(\Sigma).\]
\end{prop}

The real enumerative invariants of  toric surfaces can also be computed via tropical  geometry \cite[Thm. 6]{Mi}. 
For any real line bundle $L$ and any $\delta\ge 0$, the tropical Welschinger 
invariant $W^{L,\delta}_{\text{trop}}$ counts
simple tropical curves in $C$ in $|L|$ passing through $\dim|L|-\delta$  points in 
$\R^2$ in tropically general position. 
Here the tropical curves $C$ are counted with the {\it Welschinger multiplicity $r(C)$}: 
\[
r(C) = \!\!\!\! \prod_{\mathrm{vertices}\, v}\!\!\!\! r(v), \;\;\;\;\;\;\;\;\;\;\;\;\; r(v) = \begin{cases} (-1)^{(m(v)-1)/2}& m(v) \hbox{ odd},\\
0& m(v) \hbox{ even}.\end{cases}
\]

The irreducible Welschinger invariants $W^{L,\delta}_{0,\text{trop}}$ are defined by summing over only
irreducible curves. It is proven in \cite{IKS} that this is  independent of the points as long as they are in tropical general position.  Finally \cite{Mi} shows that there exists a set $\Sigma$ of $\dim |L|-\delta$ real points of $S$, so that 
$W^{L,\delta}(\Sigma)=W^{L,\delta}_{\text{trop}}$, and $W^{L,\delta}_0(\Sigma)=W^{L,\delta}_{0,\text{trop}}$.

If $S$ is $\P^2$ or a rational ruled surface, there is a recursion for the tropical Welschinger invariants \cite{IKS}. 
We write it in a modified form which makes the close relation to the recursion for the Severi degrees more 
evident.\footnote{
Also the multiplicity assigned in \cite{IKS} differs from those given above which we have taken from \cite{Mi}, 
but it can be shown they are equivalent.
}

\begin{defn} A sequence $\alpha=(\alpha_1,\alpha_2,\ldots)$ is called {\it odd} if 
$\alpha_i=0$ for all even $i$.
\end{defn}

\begin{rec}
Let  $L$ be a line bundle on $S$ and let $\alpha$, $\beta$ be odd sequences with $I\alpha+I\beta=EL$, and a let $\delta\ge 0$ be an integer.
With the same notations and assumptions and initial values as for the relative Severi degrees in \recref{CHrec}
the relative tropical Welschinger invariants $W^{L,\delta}_{\text{trop}}(\alpha,\beta)(y)$ are given by the following recursion formula:
if $\gamma(L, \beta,\delta)>0$, 
\begin{equation}\label{Wrec}
\begin{split}
W^{L,\delta}_{\text{trop}}&(\alpha,\beta)=\sum_{k\ \text{odd} :\beta_k>0} (-1)^{(k-1)/2} \cdot W^{L,\delta}_{\text{trop}}(\alpha+e_k,\beta-e_k)(y)
\\&+\sum_{\alpha',\beta',\delta'}
\prod_{i\ \text{odd}} \big((-1)^{(i-1)/2}\big)^{\beta_i'-\beta_i}\binom{\alpha}{\alpha'}\binom{\beta'}{\beta} W^{L-E,\delta'}_\trop(\alpha',\beta').
\end{split}
\end{equation}
Here the second sum runs through all  odd sequences $\alpha',\beta'$ and all $\delta'$ satisfying \eqref{relcong}.
\end{rec}
We put $W^{L,\delta}_{\text{trop}}:=W^{L,\delta}_{\trop}(0,LE)$, and  in the case $S=\P^2$, $W^{d,\delta}_\trop(\alpha,\beta)=W^{dH,\delta}_\trop(\alpha,\beta)$, $W^{d,\delta}_\trop=W^{dH,\delta}_\trop(0,d)$.

\vspace{4mm}

Note the following specialization:

\begin{equation} [k]_{-1} = \frac{y^{k/2}-y^{-k/2}}{y^{1/2}-y^{-1/2}}\Bigm|_{y=-1}
=\begin{cases} (-1)^{(k-1)/2}& k \hbox{ odd},\\
0& k\hbox{ even}.\end{cases}
\end{equation}

In particular, the recursion for the refined Severi degrees interpolates between the
Caporaso-Harris recursion for Severi degrees and the Itenberg-Kharlamov-Shustin recursion
for tropical Welschinger invariants.  Thus, 

\begin{prop} \label{prop:specialization}
$\overline{N}^{L,\delta}(\alpha, \beta)(1) = n^{L,\delta}(\alpha,\beta)$ and
$\overline{N}^{L,\delta}(\alpha, \beta)(-1) = W^{L,\delta}_{\trop}(\alpha,\beta)$.
\end{prop}

In \cite{BG}, Block and the first author relate the refined Severi degrees to tropical geometry and study them by the methods of tropical geometry.  They introduce
the refined multiplicity $M(v):=[m(v)]_y$, which specializes to $m(v)$ at $y=1$ and to $r(v)$ at $y = -1$.  
Then the refined tropical Severi degrees $\overline N^{L,\delta}_\trop(\Sigma)$ are defined
by counting curves with multiplicity $M(C) = \prod M(v)$.   Note this definition applies to any smooth toric surface. 
It is shown that,
for $S=\P^2$ or a rational ruled surface, and $\Sigma$ a ``vertically stretched'' configuration of points, 
the $\overline N^{L,\delta}_\trop(\Sigma)$ satisfy the recursion \eqref{norrefrec}.  Thus
$\overline N^{L,\delta}=\overline N^{L,\delta}_\trop(\Sigma)$.  

This is the analogue of the tropical proof of the Caporaso-Harris recursion formula in \cite{GM}, and like the original proof of Caporaso-Harris it can be viewed as a proof by degeneration. For  a vertically stretched configuration of points the tropical curve degenerates, so that it can be described in terms of tropical curves of lower degree, and this gives the recursion both for the Severi degrees and the refined Severi degrees.


Itenberg and Mikhalkin have in the meantime shown in \cite{IM}  that $\overline N^{L,\delta}_\trop(\Sigma)$ 
is independent of $\Sigma$, and so we drop it from the notation. 
For $S = \P^2$ or a rational ruled surface, Conjecture \ref{ref-sev} then implies that
the $\overline N^{L,\delta}_{\trop}$ agree with the $\overline N^\delta_{\delta,[S,L]}$ when 
$L$ is $\delta$ very ample.  More generally one expects:

\begin{conj} \label{conj:BG} \cite{BG}
Let $S$ be a smooth projective toric surface and $L$ a real line bundle on $S$. 
If $L$ is $\delta$-very ample, then $\overline N^{\delta}_{\delta,[S,L]} = \overline{N}^{\delta, L}_{\trop}$.
\end{conj}

Using the refined multiplicity, in \cite{BG} the  $\overline N^{L,\delta}_\trop(\Sigma)$ are studied using methods similar to those employed in \cite{Blo} for the nonrefined Severi degrees.  In particular it is shown that,  for $L$ sufficiently ample with respect to $\delta$, they are given by refined node polynomials, and the Conjecture says that these agree with the $\overline N^{\delta}_{\delta,[S,L]}$.

According to Conjecture \ref{conj:BG} and Proposition \ref{prop:specialization}, we expect:
\begin{conj}
	Let $L$ be a $\delta$ very ample real line bundle on a toric surface $S$.  Then 
 $\overline{N}^{\delta}_{\delta,[S,L]}(-1) = W^{\delta, L}_{\trop}$. 
 \end{conj}

For convenience we record the corresponding specialization of Conjecture \ref{Gconj} at $y = -1$. 
 Consider
$\eta(\tau):=q^{1/24}\prod_{n>0}(1-q^n)$  the Dirichlet eta function and write $$
\overline G_2(\tau):=G_2(\tau)-G_2(2\tau)
=\sum_{n>0} \left(\sum_{d|n, \ d \text{ odd}} \frac{n}{d} \right)q^n.$$

\begin{conj}
\label{GSconj}
\begin{equation}
\label{GSform}
\sum_{\delta\ge 0}
\overline N^\delta_{\delta, [S,L]}(-1)\overline G_2(\tau)^\delta=\frac{( \overline G_2(\tau)/q)^{\chi(L)}B_1(-1,q)^{K_S^2}B_2(-1,q)^{LK_S}}
{\big(\eta(\tau)^{16}\eta(2\tau)^4\, D\overline G_2(\tau)/q^2\big)^{\chi(\oo_S)/2}}.
\end{equation}
\end{conj}
This conjecture has been checked modulo $q^{15}$ and the coefficients of $B_1(-1,q),B_2(-1,q)$ have been determined modulo $q^{15}$. (These computations are numerically easier than those involving the indeterminate $y$, thus we get to a higher order in $q$).  The values of the series $B_i$ are computed to be:
\begin{align*}
B_1(-1,q)&=1 - q- q^2 - q^3 + 3q^4 + q^5 - 22q^6 + 67q^7 - 42q^8 - 319q^9\\& + 1207q^{10}- 1409q^{11} - 3916q^{12} + 20871q^{13} - 34984q^{14}+O(q^{15}),\\
B_2(-1,q)&=1 + q + 2q^2 - q^3 + 4q^4 + 2q^5 - 11q^6 + 24q^7 + 4q^8 - 122q^9\\& + 313q^{10} - 162q^{11} - 1314q^{12}+ 4532q^{13} - 4746q^{14}+O(q^{15}). \\
\end{align*}

When $S = \proj^2$ or a rational ruled surface, the Severi degrees $n^{\delta, L}$
agree with the universal numbers $n^{\delta}_{\delta,[S,L]}$ somewhat beyond the regime where $L$ 
is $\delta$ very ample.  Specifically, it is conjectured in \cite{G} and proven in \cite{KS} that it suffices for 
the general $\proj^\delta \subset |L|$ to contain no nonreduced curves, and no curves containing
components with negative self intersection.  We expect the same to hold for the comparison between
refined Severi degrees $N^{\delta, L}$ and the universal numbers $N^{\delta}_{\delta,[S,L]}$, and 
a fortiori for the specialization at $-1$.  However for this specialization more seems to be true:

\begin{conj}\label{welam}
Assume $S=\P^2$ or $S=\Sigma_e$, and the following
subloci of $|L|$ have codimension more than $\delta$: 
(1) the nonreduced curves with a component of multiplicity at least $3$, 
(2) curves containing a component with negative self intersection. Then 
$W^{L,\delta}_\trop=\overline N_{\delta,[S,L]}^\delta(-1)$.  Explicitly the condition amounts
to:
\begin{enumerate}
\item On $\P^2$ we have $W^{d,\delta}_\trop=\overline N^\delta_{dH}(-1)$ if $\delta\le 3d-3$.
\item On $\P^1\times \P^1$ we have $W^{nF+MG,\delta}_\trop=\overline N_{nF+mG}^\delta(-1)$
if $\delta\le 3\min(n,m)$.
\item On $\Sigma_e$ with  $e>0$, we have
$W^{nF+mE,\delta}_\trop=\overline N_{nF+mE}^\delta(-1)$, if 
$\delta\le \min(3m,n-em)$.
\end{enumerate}
\end{conj}

Using the recursion formula \eqref{Wrec} this conjecture has been checked for $d,\delta\le 14$.
Assuming (1) of \conjref{welam}, and using the recursion formula \eqref{Wrec}, the conjecture
\ref{GSconj} has been checked modulo $q^{67}$  and $B_1(-1,q)$ and $B_2(-1,q)$ have been determined modulo $q^{67}$. Note that the recursion for $W^{d,\delta}_\trop(\alpha,\beta)$ is much more efficient than those of the $N^{d,\delta}(\alpha,\beta)(y)$ or $n^{d,\delta}(\alpha,\beta)$ because only odd 
sequences $\alpha$ and $\beta$ occur.

\vspace{4mm}

We have seen $W^{d,\delta}_\trop = 
\overline{N}^{d,\delta}(-1)$.  In the sufficiently ample setting,  taking a linear system $\proj^\delta$
determined by subtropical point conditions, and assuming all conjectures,
this implies
$N^\delta_{\cC/\proj^\delta}(-1) =  (-1)^\delta  \overline{N}^\delta_{\cC/\proj^\delta}(-1) = 
 n^{\delta,\R}_{\cC/\proj^\delta}$.  More generally we conjecture: 

\begin{conj} \label{conj:harnack}
	Let $L$ be a sufficiently ample real line bundle on a real toric surface, and let $\proj^\delta \subset |L|$ 
	be determined by a subtropical collection of point conditions.  Then the signatures
	of  the relative Hilbert schemes agree with the compactly
	supported Euler characteristics of their real loci.  That is, 
	$\chi_{1}(\cC^{[n]}_{\proj^\delta}) = \chi_c(\cC^{[n]}_{\proj^\delta}(\R))$. 
	Making the BPS change
	of variables, it follows that $N^\delta_{\cC/\proj^\delta}(-1) = n^{\delta,\R}_{\cC/\proj^\delta}$.
\end{conj}

More generally one may consider the question: 

\begin{quest}
Let $X$ be a smooth real variety.  When is
$\chi_{1}(X) = \chi_c(X(\R))$? 
\end{quest} 

This has been the subject of some classical study, one general result being that
for a ``M-variety'', i.e. one for which the total dimension of the $\Z/2\Z$ cohomology is
equal for the real and complex locus, the equality holds modulo 16 \cite{DK}.  

Evidently the equality holds for any variety whose class in the Grothendieck ring of varieties over $\R$ lies inside
$\Z[\mathbb{A}^1_{\R}]$.  In particular, $\mathbb{RP}^n$, toric surfaces, 
and Hilbert schemes of points on toric surfaces qualify.  The relative Hilbert schemes are cut
out of a product of these by a section of a vector bundle, and the signature behaves predictably under
taking such sections.  Thus, to study Conjecture \ref{conj:harnack} it would also suffice 
to give criteria for the Euler characteristic of the real locus to exhibit the same behavior. 
In \cite{BB} it is shown this holds in a tropical sense for complete intersections in toric varieties.

\begin{rem}
The conjectured relation between the refined invariants at $y=-1$ and the Welschinger invariants is in a sense
a global analogue of a conjecture of van Straten \cite[Conj. 4.6, 4.7]{vS} (see also Theorems
\ref{duco} and \ref{thm:realmult} above and
the nearby discussion). 
\end{rem}

\end{document}